\documentclass{ieeeaccess}
\usepackage{xcolor} 
\usepackage[T1]{fontenc}
\usepackage[utf8]{inputenc}
\usepackage{cite}
\usepackage[hyphens]{url}

\usepackage[hidelinks,bookmarks=false]{hyperref}
\usepackage{amsmath,amssymb,amsfonts}
\usepackage{algorithmic}
\usepackage{graphicx}
\usepackage{textcomp}
\usepackage[font=small,labelsep=colon]{caption}
\usepackage[vlined,ruled]{algorithm2e}
\usepackage{array}
\newcolumntype{?}{!{\vrule width 1.2pt}}
\usepackage{diagbox}

\usepackage{stfloats}
\usepackage{comment}
\usepackage{longtable}
\usepackage{booktabs}
\usepackage{multirow,multicol,makecell}
\usepackage{tabularx}
\usepackage{enumitem}

\def\BibTeX{{\rm B\kern-.05em{\sc i\kern-.025em b}\kern-.08em
    T\kern-.1667em\lower.7ex\hbox{E}\kern-.125emX}}

\begin{document}
\history{}
\doi{}

\title{Evaluating DV/CV-QKD Architectures for SAFE Long-Term Secure Storage: A Risk Model and
ILP-Based Cost Optimization Approach}
\author{\uppercase{A. TASDIGHI}\authorrefmark{1} and 
\uppercase{R. Alléaume}\authorrefmark{1,2}}
\address[1]{Télécom Paris, LTCI, Institut Polytechnique de Paris, Inria, 19 Place Marguerite Perey, 91120 Palaiseau, France (email: alireza.tasdighi@telecom-paris.fr, romain.alleaume@telecom-paris.fr)}
\address[2]{ Inria Saclay, 1 Rue Honoré d’Estienne d’Orves, 91120 Palaiseau, France (email: romain.alleaume@telecom-paristech.fr)}
\tfootnote{This article extends the preliminary version presented at the 2026 International Conference on Quantum Communications, Networking, and Computing (QCNC) \cite{Tasdighi2026}.}

\markboth
{Author \headeretal: Preparation of Papers for IEEE Transactions on Quantum Engineering}
{Author \headeretal: Preparation of Papers for IEEE Transactions on Quantum Engineering}

\corresp{Corresponding author: A. TASDIGHI (email: alireza.tasdighi@telecom-paris.fr).}

\begin{abstract}
This paper presents a unified cost-modeling and optimization framework designed to evaluate hybrid discrete-variable/continuous-variable quantum key distribution (DV/CV-QKD) infrastructures operating within the SAFE (Secure and Efficient) long-term storage (LTSS) protocol. Our methodology jointly addresses the information-theoretic and computational dimensions of cryptographic durability over multi-decade horizons by coupling a quantitative risk model with a stochastic integer linear programming (ILP) formulation. Going beyond idealized physics-informed limits, the framework incorporates realistic next-generation industrial implementations and multiplexed coexistence specifications, together with a sample-average approximation (SAA) pipeline, to determine cost-optimal and structurally feasible SAFE topologies under transmission-budget and security constraints. The proposed framework is evaluated on both randomized synthetic deployments and realistic metropolitan-scale QKD infrastructures, including the Paris Metro-Scale and Greater Paris networks, to characterize the interplay between network topology, QKD modality, and deployment cost while enforcing a target global compromise tolerance $\epsilon$. Numerical results show that the economically optimal architecture is highly scenario-dependent and that both hybrid and homogeneous modality allocations naturally emerge from the optimization. More importantly, the analysis reveals a non-monotonic relationship between the minimum required QKD link capacity and the infrastructure cost: contrary to intuition, increasing the minimum required per-link secret-key rate may reduce the overall deployment cost by enabling a global transition toward more economical CV-based solutions, whereas further capacity increases do not necessarily yield additional savings. These findings pave the way towards cost-effective long-term secure data storage implemented with modern commercially available QKD technologies.
\end{abstract}


\titlepgskip=-15pt

\maketitle

\section{Introduction}

In an era where digital services must guarantee confidentiality over decades, long-term secure storage (LTSS) is a growing industrial priority. Applications such as archival medical records, intellectual property vaults, and national digital archives demand systems whose encryption cannot be retroactively broken even if future advances render current cryptographic assumptions obsolete. Recent deployments of metropolitan-scale quantum networks underscore this trend: for instance, Geneva inaugurated its first 262 km quantum fiber network GQN\cite{GQN2025} to support quantum communications and cryptographic experiments in situ, while in Germany the QuNet\cite{QuNET2025} initiative links fiber sites across Berlin as a testbed for quantum network integration. In France also, the ParisRegionQCI\cite{OrangeQCI2023} backbone interconnects four metropolitan quantum nodes across an 80 km fiber span with two trusted relay hops, over which both continuous-variable (CV) and discrete-variable (DV) quantum key distribution (QKD) systems have been deployed. These examples suggest both growing technical feasibility and institutional momentum for quantum-secure infrastructures in the real world.

From a cryptographic perspective, classical authenticated key exchange (AKE) (e.g., RSA\cite{Rivest1978RSA} or symmetric AKE with AES\cite{NIST2001AES}) offers only \emph{computational security} (CS): if the underlying hard problems are later broken, past and future sessions may be exposed (no forward or post-compromise resilience) ~\cite{mosca2013quantum,alleaume2021HDR}. By contrast, QKD provides \emph{information-theoretic security} (ITS), decoupling confidentiality from computational assumptions. When QKD-derived keys are combined with one-time pad (OTP) encryption, possibly using computationally-secure authentication the resulting scheme attains \emph{everlasting confidentiality}, i.e. security against a computationally unbounded adversary after protocol execution. On the other hand, combinations such as QKD$+$AES, do not significantly strengthen security: any computational vulnerability (like quantum computational power) in data encryption can still compromise the stored data chain\cite{Bonnetain2019}. Thus, for metropolitan-scale LTSS, QKD$+$OTP remains a uniquely valid option — provided the system is engineered for efficiency and cost-risk balance.

A key challenge is enabling scalable and cost-effective LTSS architectures that deliver ITS confidentiality without unmanageable operational overhead. Several protocols have been proposed, including LINCOS (with COPRIS)\cite{Bouchmann2017LINCOS}, PROPYLA\cite{Geihs2018PROPYLA}, ELSA / MCELSA\cite{Bouchmann2020ELSA}, SPSS+TTP\cite{Fujiwara2021LTSS}, MULTISS\cite{Prevoste2024MUTLISS}, and SAFE\cite{Bouchmann2020SAFE}. The designs differ in how they separate (or unify) confidentiality and integrity, manage renewable commitments (such as Pedersen’s Commitment Scheme (PCS) \cite{Pedersen1991}), deploy trusted components (Trusted Execution Environment (TEE) or Trusted Third Party (TTP)), and scale to multiple nodes. We compress the detailed survey to a single paragraph in Section~\ref{subsec:sota}, noting that SAFE stands out with a TEE-centric, star-like topology that unifies integrity and confidentiality within a single enclave, thereby minimizing the number of ITS channels and simplifying trust management.

In this work, we propose a risk-aware optimization framework for hybrid DV/CV-QKD architectures under the SAFE protocol. We develop a joint compromise-risk model and a scenario-based integer linear programming (ILP) solver that decides which QKD modalities and network topologies minimize cost while satisfying a global security tolerance \(\epsilon\). Our methodology flexibly adapts to physical key-rate models, cost uncertainties, and renew cadence. Through simulations over random metropolitan layouts, we demonstrate that mixed DV/CV configurations, especially when combined with multiplexing, yield significant cost savings compared to single-modality deployments, while preserving ITS-level confidentiality.

Compared with the preliminary conference version \cite{Tasdighi2026}, this article expands the SKR modeling of QKD modalities in Section \ref{sec:QKD-mods} through additional industrial specifications, coexistence-aware models, and comparative figures/tables. The minimum per-link secret-key-rate requirement is further formalized and analyzed in Section \ref{sec:SKR_requirement}, while additional metropolitan-scale benchmark studies using realistic industrial QKD imperfections and the proposed ILP framework are provided through extended numerical evaluations. The remainder of this paper is organized as follows: Section~\ref{sec:prelim} presents the necessary cryptographic and QKD preliminaries. Section~\ref{sec:risk} develops the compromise-risk framework, and Section~\ref{sec:ilp} formalizes the optimization problem. Section~\ref{sec: integ} describes the integrated risk–optimization pipeline. Section~\ref{sec:synth} synthesizes the simulation results. We conclude in Section~\ref{sec:concl}.

\section{Preliminaries and backgrounds}
\label{sec:prelim}

This section introduces the security notions, cryptography primitives and brief background used throughout the paper.

\subsection{Security notions: ITS vs. CS}
We distinguish two standard notions of security: 1) ITS, that guarantees certain secrets remain statistically (or perfectly) hidden regardless of the adversary's computational power. Examples: OTP, secret sharing with unconditional privacy. 2) Computationally Secure, where security depends on assumed hardness of computational problems (e.g., RSA). Guarantees hold only against bounded adversaries within feasible time/resources.
\subsection{Shamir's Secret Sharing (SSS)~\cite{Shamir1979SSS}}
Let $\mathbb{F}_q$ be a finite field with $q = 2^r$ and $q > m$. 
The original document $M \in \{0,1\}^{|M|}$ is split into $N = \lceil |M|/r \rceil$ chunks of size $r$ bits, where chunk $\ell$ $(1 \leq \ell \leq N)$ is represented as an element of $\mathbb{F}_q$. 
To distribute chunk $\ell$ among $m$ storage nodes with reconstruction threshold $t$, the following procedure is applied:

\textbf{(1)} Choose random coefficients $a_1, \dots, a_{t-1} \in \mathbb{F}_q$ and define the polynomial
$f_\ell(x) = \sigma_{0,\ell} + a_1 x + \cdots + a_{t-1} x^{t-1}$,
  where $\sigma_{0,\ell} \in \mathbb{F}_q$ is the chunk value.
  
\textbf{(2)} For each storage node $\text{S}_i$ (with distinct identifier $i \in \mathbb{F}_q^*$), compute the share $\sigma_{i,\ell} = f_\ell(i)$.
The full share held by node $\text{S}_i$ is 
$\sigma_i = (\sigma_{i,1} \,\Vert\, \sigma_{i,2} \,\Vert\, \cdots \,\Vert\, \sigma_{i,N}) \in (\mathbb{F}_q)^N$.

Any $t$ shares $\{(i, \sigma_{i,\ell})\}$ are sufficient to reconstruct $\sigma_{0,\ell} = f_\ell(0)$ using Lagrange interpolation, while any set of fewer than $t$ shares reveals no information about $m_\ell$, thereby ensuring ITS. 
Although in principle $2 \leq t \leq m$, in practice many LTSS systems assume $t = \lceil (m+1)/2 \rceil$, known as the \textit{honest majority} assumption, which balances confidentiality (resisting up to $t-1$ compromises) and availability (tolerating up to $m-t$ node failures).

\subsection{Proactive Secret Sharing (PSS)}
PSS~\cite{Herzberg1995PSS} is to defend against mobile adversaries who may compromise different nodes over time, so that the shares are periodically refreshed without altering the underlying secret. 
In each refresh round (occurring every $\Delta^{\text{Refresh}}$ years), random degree-$(t-1)$ polynomials $g_\ell(x)$ with $g_\ell(0)=0$ are sampled. 
Each node $\text{S}_i$ receives an update value $p_{i,\ell}=g_\ell(i)$, and updates its share as
\[
\sigma_{i,\ell} \leftarrow \sigma_{i,\ell} + p_{i,\ell}.
\]
This ensures that the secret chunks $m_\ell$ and threshold $t$ remain unchanged, but previously compromised shares become statistically independent of the new ones. 
Thus, confidentiality is preserved as long as at most $t-1$ nodes are corrupted between two successive refreshes.

\subsection{QKD Modalities $\&$ Their Key Rates}\label{sec:QKD-mods}
QKD enables two remote parties to establish ITS keys over optical fiber. Each system employs two parallel links: a quantum optical channel for transmitting quantum states and a classical authenticated channel for post-processing tasks such as basis sifting, error correction, and privacy amplification. 
The efficiency of this classical post-processing, quantified by the \textit{reconciliation efficiency} ``$\beta$" directly impacts the achievable SKR.

Two main QKD modalities are considered. DV-QKD, exemplified by the BB84 protocol\cite{BB84}. It supports long-distance transmission (beyond $200$~km on a single fiber span) but requires sensitive single-photon detectors. On the other hand, CV-QKD is well suited for deployment in coexistence with classical traffic, multiplexed over the same fiber\cite{kumar2015coexistence}. It employs homodyne detection, achieves higher raw key rates at metropolitan scales, but is more sensitive to loss. 

When combined with multiplexing (MUX) techniques like wavelength division multiplexing (WDM), both modalities suffer from crosstalk and excess noise. DV$+$MUX systems experience moderate degradation, whereas CV$+$MUX configurations are more noise sensitive, leading to significant reductions in both SKR and the secure distance.

Transitioning from fundamental physical limits to phenomenological implementations that represent current and next-generation industrial deployments, three-tiers for SKRs could be considered:

\subsubsection{Theoretical Upper Bound}
\;\;The fundamental secret capacity limit of a repeaterless lossy channel is established by the Pirandola-Laurenza-Ottaviani-Banchi (PLOB) bound \cite{Pirandola:2017kzd}. The PLOB bound defines the maximum achievable key rate over a purely lossy channel as a function of the link transmission distance $d$ (in km) and is expressed as follows:
\begin{equation}
\label{InfTheoSKR}
\mathrm{SKR}_{\mathrm{PLOB}}(d) = C_{\mathrm{norm}} \cdot \left[ -\log_2 \left( 1 - 10^{-\alpha_\mathsf{eff}\;d} \right) \right],
\end{equation}
where $\alpha_\mathsf{eff} = 0.02$ characterizes a standard optical fiber attenuation of $0.2~\mathrm{dB/km}$. The scaling coefficient $C_{\mathrm{norm}}$ is a normalization constant introduced to map the channel-capacity bound into operational throughput units of bits per second ($\mathrm{bits/s}$), providing an idealized baseline reference for network evaluation.

\subsubsection{Trivial Baseline Model}
\;\;The SKR–distance model used in this tier is physics-informed but trivial yet, intended to capture relative performance trends rather than exact laboratory rates. It combines exponential attenuation with distance (modeling fiber loss and device inefficiency) and a modality-specific SKR at zero distance $\text{SKR}_0^{(\cdot)}$, where $\cdot\in \lbrace \text{CV}, \text{DV}, \text{CV}+\text{MUX}, \text{DV}+\text{MUX}\rbrace$. Mathematically, for a link of length $d$, this ideal experimental baseline is modeled as:
\begin{equation}
\label{eq:capacity}
\text{SKR}^{(\cdot)}(d) = \, \text{SKR}_0^{(\cdot)}\,\big[10^{-\alpha_{\mathrm{eff}}\, d} - 10^{-\alpha_{\mathrm{eff}}\, d^{(\cdot)}_{0}}\big]_+ ,
\end{equation}
where $[\cdot]_+ = \max(\cdot, 0)$, and $d^{(\cdot)}_{0}$ denotes the operational cut-off distance beyond which the rate vanishes. The parameters Per-modality $\text{SKR}_0^{(\cdot)}$ and $d^{(\cdot)}_{0}$ are reported in the two left columns of Table~\ref{tab:modality_overview}. These values are consistent in order of magnitude with recent experimental results on metropolitan links~\cite{bebrov2024relation, Pirandola:20, Pi:23, kumar2015coexistence}. The relevant SKR–distance curves for this trivial baseline are illustrated and compared with the theoretical upper bound in Fig.~\ref{fig:Ideal_Industrial_SKR} (a).

\subsubsection{Industrial Implementation $\&$ Multiplexed Coexistence}
\;\;This tier models performance degradation caused by hardware integrated component losses, imperfect error reconciliation, and ambient dark count or standard optical noise floors as below:
\begin{equation}
\label{IndustSKR}
\text{SKR}^{(\cdot)}(d) =\left[ \text{SKR}^{(\cdot)}_0 \cdot 10^{-\alpha_\mathrm{eff}\; d} \right] - \left( \beta^{(\cdot)} \cdot \mathbf{I}_\text{noise} \right),
\end{equation}
where $\mathbf{I}_\text{noise}$ defines the empirical noise-floor penalty in $\mathrm{bits/s}$ and fiber attenuation $\alpha_\mathrm{eff}$ increases during the coexistence of the classical MUX. 

\begin{figure*}[htb!]
\centering
\begin{tabular}{c@{~}c}
\includegraphics[width=8.1cm, height=6.5cm]{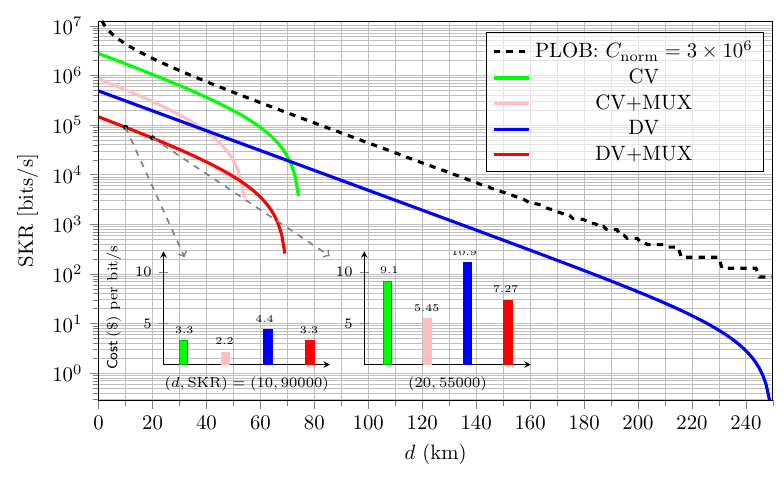} &
\includegraphics[width=8.19cm, height=6.41cm]{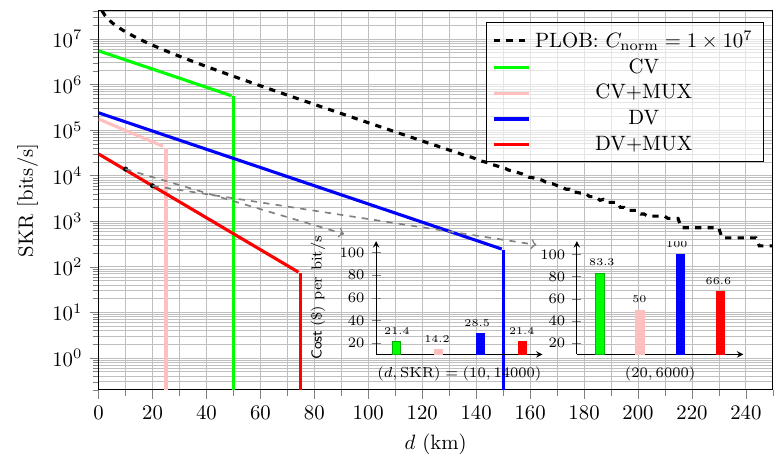} \\ 
~~(a) & ~~(b)
\end{tabular}
\caption{SKR as a function of fiber distance $d$ alongside the PLOB bound under (a) baseline setup parameters and (b) industrial specifications. Insets: Comparative infrastructure deployment cost per bit/s evaluated under the 10-year depreciation model presented in \ref{QKD_deploy_cost_model}. The black markers designate the fixed operational coordinate points $(d, \text{SKR})$ selected to benchmark the economic efficiency of all four modalities prior to their respective distance cut-off thresholds.}
\label{fig:Ideal_Industrial_SKR}
\end{figure*}

The underlying technical parameters mapping onto modern commercial equipment and industrial pilot deployments are consolidated in Table~\ref{tab:modality_overview}, with their corresponding performance curves illustrated in Fig.~\ref{fig:Ideal_Industrial_SKR} (b).

\begin{table*}[htb!]
\centering
\caption{Trivial and Industrial QKD Modality Parameters and Implementers}
\label{tab:modality_overview}
\begin{tabular}{?@{}c@{}|@{}c@{}?@{}c@{}?@{}c@{}|@{}c@{}|@{}c@{}|@{}c@{}|@{}c@{}|@{}c@{}?}
\noalign{\hrule height 1.2pt}
\multicolumn{2}{?@{}c@{}?}{\textbf{Trivial Specifications}}& \multirow{2}{*}{\textbf{ Modality Type }} & \multicolumn{6}{c@{}?}{\textbf{Industrial $\&$ Coexistence Specifications}}\\ \noalign{\global\arrayrulewidth=1.2pt} \cline{1-2}\cline{4-9} \noalign{\global\arrayrulewidth=0.4pt}

$\;\text{SKR}^{(\cdot)}_0$\; & \textbf{$d^{(\cdot)}_0$} &  & $\;\text{SKR}^{(\cdot)}_0\;$ & $\mathbf{\alpha}_\text{eff}$ & \textbf{$d^{(\cdot)}_0$} &  \textbf{$\;\;\;\;\:\:\:\beta^{(\cdot)}\:\:\:\;\;\;\;$} & $\;\;\;\;\;\:\:\mathbf{I}_\text{noise}\:\:\;\;\;\;\;$ &  \textbf{Implementers} \\ \hline
$\:1e5-5e5\:$ & $\:250-280\:$ & DV  & $1.8e5$ & $0.02$ & $\;150\;$ & $1.10-1.15$ & ${10}^{-8}-{10}^{-7}$ &\;Toshiba, QUKY, QuantumCTek, IDQ\; \\ \hline
$\:8e4-1.6e5\:$ & $\:60-80\:$ & DV+MUX & $2e4$ & $\;0.035\;$ & $75$ & $1.20$ & ${10}^{-5}-{10}^{-4}$ & Toshiba, IDQ \\ \hline
$\:3e6-1e7\:$ & $\:75-95\:$ & CV & $5e6$ & $0.02$ & $50$&  $0.94-0.96$ & $0.01-0.02$ & LuxQuanta, Huawei, Thales \\ \hline
$\:1e6-3e6\:$ & $\:55-75\:$& CV+MUX & $1.5e5$ & $0.025$ & $25$&  $0.80-0.85$ & $0.10-0.15$ & LuxQuanta, ADVA, Orange (Pilots) \\ \noalign{\hrule height 1.2pt}
\end{tabular}
\end{table*}
\vspace*{0.7cm}

\subsection{Per QKD Link Deployment Cost Model}\label{QKD_deploy_cost_model}
The ten-year hardware depreciation cost estimate  per-link is roughly defined as affine distance functions $C_{i,j}^{(\cdot)}$ $=C^{cl}d_{i,j}+C^{(\cdot)}$ with fixed offsets of $C^\text{DV}$ $= C^\text{DV$+$MUX}$ $= 200\ \text{k}\$$, $C^\text{CV}$ $= C^\text{CV$+$MUX} = 100\ \text{k}\$$ where, $d_{i,j}$ is the distance (in km) between the nodes $i$ and $j$. One-year rental cost (per km) of classical fiber link is also estimated by $C^\text{cl} = 1\ \text{k}\$$. Perturbed by multiplicative interval factor $\gamma=[\underline{\gamma},\overline{\gamma}]$ we define $[\underline{C^{(\cdot)}_{i,j}}, \overline{C^{(\cdot)}_{i,j}}]$ $= C^{(\cdot)}_{i,j}\ \gamma$. This interval is used for per-modality-cost random sampling (see Sec.~\ref{sec: integ}, Alg.~\ref{alg:integrated}).

It is worth noting that the considered cost function here serves as trivial physics-based input to explore cost, feasibility, and modality trade-offs under stochastic deployment conditions. Moreover, the architectural cost model formalized above is benchmarked and illustrated for both baseline and industrial QKD-link deployment contexts via the inset plots in Fig. 1. Specifically, to evaluate the practical economic impact of these expenditure formulations, we map the 10-year depreciation cost $C_{i,j}^{(\cdot)}$ against the achieved $\text{SKR}^{(\cdot)}(d)$). This yields a normalized cost-per-performance metric, expressed as the infrastructure deployment cost per unit of throughput ($\$ \text{ per bit/s}$). By anchoring this cost model to the physical layer simulations, the inset figures provide a direct economic comparison of the four QKD modalities under identical operational constraints.

\subsection{SAFE Protocol}
\label{subsec:SAFE}
SAFE provides a TEE-coordinated framework for LTSS, jointly maintaining confidentiality and integrity through PSS and computational authentication. A single provisioned TEE orchestrates the protocol, while $m$ storage providers $\{\text{S}_i\}_{i=1}^m$ hold encrypted shares. The design is deliberately \emph{TEE-agnostic}: the enclave serves as a replaceable coordination point rather than a long-term data holder. This enables secure \emph{TEE migration} across hardware life cycles—when a platform becomes obsolete or untrusted, SAFE’s $\mathsf{Renew}$ phase transfers authority to a newly attested TEE without revealing stored content.

\paragraph{$\mathsf{Share}$, $\mathsf{Renew}$, and $\mathsf{Reconstruct}$}
SAFE employs PSS with periodic resharing of data fragments. The secret document is divided into $N$ chunks ($N\in\mathbb{N}$), each node $\text{S}_i$ storing $\sigma_i\in(\mathbb{F}_q)^N$. Under TEE supervision, three core operations are executed:

\textbf{(1) $\mathsf{Share}$:} The TEE generates random degree-$(t\!-\!1)$ polynomials to compute $\{\sigma_i\}_{i=1}^m$ from the secret $M\!\in\!\{0,1\}^{|M|}$ with identifier $\text{ID}_M$. For each $\text{S}_i$, it issues an RSA-based signature $\sigma_i^{\text{RSA}}=\mathsf{Sign}_{\text{RSA}}(\sigma_i\|\text{ID}_M)$ and transmits $(\sigma_i,\sigma_i^{\text{RSA}})$ over authenticated ITS channels. Only a minimal provisioning ledger (signatures and timestamps) is retained inside the TEE, preventing exposure of raw shares.

\textbf{(2) $\mathsf{Renew}$ ($\mathsf{Reshare}$):} At each refresh interval $\Delta^{\text{Refresh}}$, the TEE samples fresh random polynomials and computes additive updates $p_{i,\ell}$ per chunk $\ell$ so that $\sigma_{i,\ell}\!\leftarrow\!\sigma_{i,\ell}+p_{i,\ell}$ for $i=1,\dots,m$. This re-randomization renders successive share sets statistically independent. During Renew, the TEE may also rotate attestation keys and reissue signatures, supporting enclave replacement while maintaining confidentiality.

\textbf{(3) $\mathsf{Reconstruct}$:} When authorized access is requested, a quorum of at least $t$ nodes transmits signed shares to a verified TEE endpoint. The TEE authenticates each $\sigma_i^{\text{RSA}}$, checks consistency, and reconstructs the original chunks via Lagrange interpolation. Reconstruction occurs exclusively within the attested environment; node-resident shares never leave their domains (see~\cite{Bouchmann2020SAFE}, Sec.~3.2).

\begin{figure}[ht]
\centering
\includegraphics[height=0.25\textheight,width=0.42\textwidth]{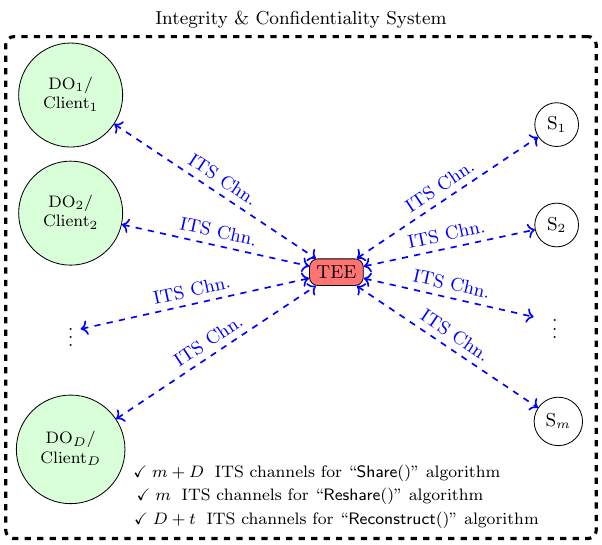}
\caption{SAFE topology with a central TEE, $D$ data owners/clients, and $m$ storage nodes. Integrity and confidentiality mechanisms are unified.}
\label{fig:SAFE_Network_Topologies}
\end{figure}

\paragraph{Network (star-like / TEE-centric)}
SAFE assumes a star topology (Fig.~\ref{fig:SAFE_Network_Topologies}), where each $\text{S}_i$ maintains a secure authenticated channel with the central TEE. The main requirements are: (i) authenticated, replay-protected communication on TEE–$\text{S}_i$ links; (ii) mutual attestation and public-key verification, prior to provisioning or reconstruction; and (iii) secure deletion of obsolete shares after $\mathsf{Renew}$. These assumptions confine the number of ITS links to $\mathcal{O}(m)$—significantly fewer than in fully pairwise commitment-based schemes—and simplify key management.

\paragraph{Integrity assumptions}
Integrity and authenticity rely on computational digital signatures (DS) (e.g., RSA-3072) generated within the TEE. Each share $\sigma_i$ carries a signature $\sigma_i^{\text{RSA}}$ whose size and renewal cadence are reflected in $\gamma_{\text{pad}}$ and $\Delta^{\text{SAFE}}$ (Section~\ref{sec:risk}). While signatures offer computational security, confidentiality remains information-theoretic due to proactive resharing. Periodic signature rotation during Renew ensures long-term verifiability without resorting to information-theoretic commitments.

\paragraph{Adversary model}
SAFE long-term secure storage defends against harvest-now-decrypt-later (HNDL)-type and mobile adversaries capable of sequential node compromises or traffic observation. Security holds provided fewer than $t$ nodes are corrupted within a renewal epoch. Since resharing yields independent shares across epochs, even unbounded adversaries cannot combine old and new fragments to reconstruct prior plaintexts, assuming secure deletion and authenticated provisioning. The TEE is honest-but-curious—faithfully executing the protocol while isolating plaintext data—and its trust is time-bounded via periodic attestation and migration. Although computational primitives may eventually weaken, information-theoretic secrecy of the shared data persists, ensuring forward security throughout the $\mathsf{Renew}$ life cycle.

\begin{figure}[ht]
\centering
    \includegraphics[width=0.5\textwidth]{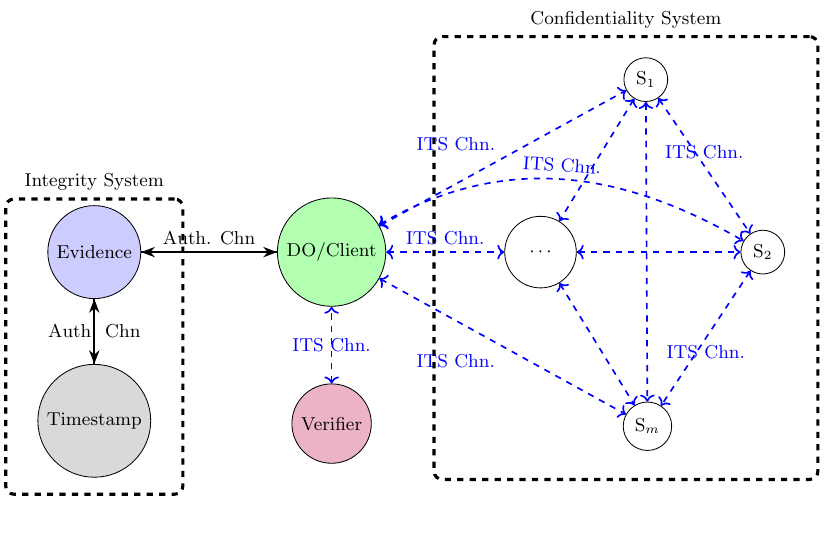} \\ 
    \caption{LINCOS $\&$ ELSA network topology with a DO/client, a verifier, an evidence service, a timestamp and, $m$ storage nodes. Integrity system is separated from confidentiality system}
    \label{fig:LINCOS_Network_Topologies}
\end{figure}

\begin{figure}[ht]
\centering   \includegraphics[width=0.48\textwidth]{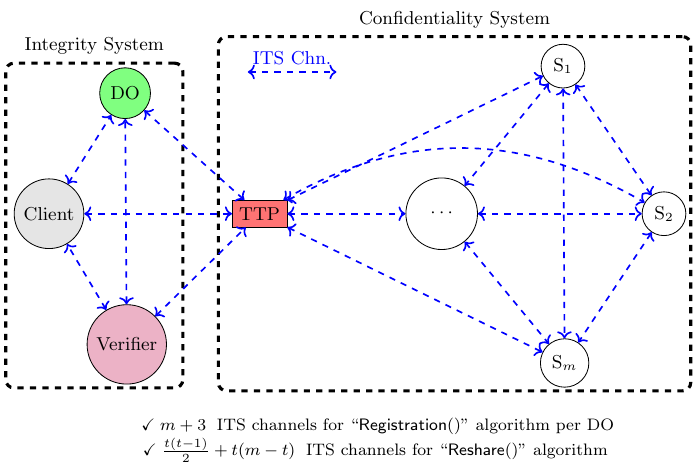}
   \caption{SPSS-TTP network topology with TTP-centric node, a DO, a client, a verifier, and, $m$ storage nodes. Integrity system is separated from confidentiality system}
    \label{fig:SPSS_TTP_Network_Topologies}
\end{figure}
\begin{figure}[ht]
\centering
  \includegraphics[width=\columnwidth,height=8.5cm]{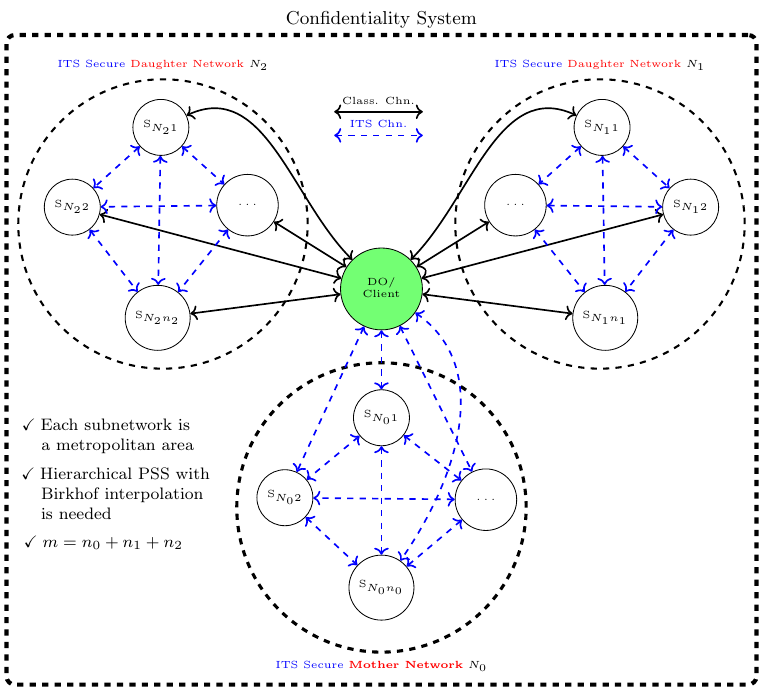} 
  \caption{MULTISS network topology with a DO/client and $m$ storage nodes. Only confidentiality system is considered}
    \label{fig:MULTISS_Network_Topologies}
\end{figure}

\subsection{SKR Requirements in SAFE Renew Procedure}\label{sec:SKR_requirement}
A minimum required link-level SKR $k_\text{min}$ must be rigorously defined within the context of the SAFE $\mathsf{Renew}$ procedure prior to formulating the risk and ILP models presented in Sections~\ref{sec:risk} and \ref{sec:ilp}. 

In SAFE, each $\mathsf{Renew}$ round consists of two sequential operational phases: 1) \textbf{Forward phase} where the TEE generates fresh secret shares and distributes them across all $m$ storage nodes 2) \textbf{Backward phase} where the TEE retrieves a qualified threshold subset of shares (at least $t$) to reconstruct and update the secret. Consequently, the total secret exchange latency incurred within a single $\mathsf{Renew}$ round is defined by:
\begin{equation}
\Delta^{\text{Refresh}} = \Delta^{\text{Forward}} + \Delta^{\text{Backward}}.
\end{equation}

Note that $\Delta^{\text{Refresh}}$ isolates the data transmission latency and excludes internal cryptographic or computational processing overheads at the TEE or storage nodes. These internal SAFE processes, detailed in Section~\ref{sec:risk} and quantified in Table~\ref{tab:SAFE_time}, execute concurrently with data transmission tasks. Because their execution times are orders of magnitude shorter than network transmission delays, they are considered negligible and are omitted from the primary latency bound.

\subsubsection{Forward Phase (Bottleneck)}
\;\;Let $k_i$ denote the achievable SKR of the link between the TEE and node $i$, and let $|M|$ denote the total size of stacks of shares. The transmission time over link $i$ is $|M|/k_i$, and thus the forward latency is governed by the slowest link:
\begin{equation}
\Delta^{\text{Forward}} = \max_{i \in \{1,\dots,m\}} \left( \frac{|M|}{k_i} \right).
\end{equation}

This expression highlights a fundamental bottleneck: the completion time of the forward phase is dictated by the \textbf{minimum required per-link SKR} in the system:
\begin{equation}
k_{\min} = \min_{i} k_i.
\end{equation}

Consequently, to satisfy a given refresh deadline $\Delta^{\text{Refresh}}$, all links must be provisioned such that
\begin{equation}
k_i \geq \frac{|M|}{\Delta^{\text{Forward}}}, \quad \forall i.
\end{equation}

\subsubsection{Backward Phase (Flexible)}
\;\; In contrast, the backward phase only requires retrieving $t$ out of the $m$ shares. Therefore the reconstruction latency can benefit from selecting the $t$ fastest links. But despite this flexibility, the overall system performance remains dominated by the forward phase and all links must satisfy the same minimum rate requirement. The slowest link imposes a hard lower bound on the achievable refresh time. In addition, infrastructure cost is governed by the weakest link.

\subsection{Other state-of-the-art LTSS systems (comparative overview)}
\label{subsec:sota}
Under the cost–risk framework, several architectural distinctions differentiate SAFE from existing LTSS solutions (illustrated in Fig.\ref{fig:SAFE_Network_Topologies}–\ref{fig:MULTISS_Network_Topologies}). First, in \emph{where ITS is realized}, COPRIS, PROPYLA, ELSA, and related PSS-based schemes achieve ITS confidentiality via secret sharing and QKD/OTP links, while integrity is enforced with renewable commitments or timestamps. SAFE attains the same ITS confidentiality but consolidates integrity proofs within a single TEE, eliminating external commitment chains. Second, in \emph{topology and trust distribution}, LINCOS, ELSA, and SPSS require many pairwise ITS links or a TTP for provisioning; MULTISS spreads risk over multiple QKD domains. In contrast, SAFE’s star-like, TEE-centric architecture confines trusted functionality to one attested enclave, reducing the needed ITS channel count from quadratic to linear in the number of storage nodes (see \cite{Bouchmann2020SAFE}, Table 1). Third, on \emph{performance and scalability}, ELSA reduces per-item timestamp overhead, PROPYLA incurs delay to mask access patterns, and MULTISS increases coordination complexity across networks. SAFE achieves comparable proactive security with far lower coordination cost, since renewal and integrity operations are centrally orchestrated in the enclave. Indeed, SAFE benefits from periodic TEE migration rather than full multi-domain reconfiguration. Finally, in \emph{integrity handling}, prior systems treat integrity and confidentiality separately—often via external verifiers—whereas SAFE unifies both within the TEE. This integration supports renewability and limits dependence on long-lived cryptographic assumptions. SAFE thus offers a lean, scalable baseline against which hybrid QKD topologies can be benchmarked.

\noindent The next section uses these notations and assumptions to build the risk model (Section~\ref{sec:risk}) and the ILP formulation for cost optimization (Section~\ref{sec:ilp}).
\section{Risk model} \label{sec:risk}
We consider a DO who wishes to store an original binary document $M\in\{0,1\}^{|M|}$ in a long-term secure manner. The document is partitioned into $N=\lceil |M|/r\rceil$ chunks of size $r$ bits, each mapped to an element of $\mathbb{F}_q$ with $q=2^r$. Using SSS with parameters $(m,t)$, each chunk is distributed among $m$ active storage nodes $\{\text{S}_i\}_{i=1}^m$ such that any set of $t$ or more nodes can reconstruct the secret, while fewer than $t$ shares reveal no information. PSS refreshes shares every $\Delta^{\text{Refresh}}$ years, thereby ensuring that compromised old shares cannot be linked to new ones.

Each share $\sigma_i\in(\mathbb{F}_q)^N$ is authenticated by an RSA-3072 DS $\sigma_i^{\text{RSA}}=\mathsf{Sign}_{\text{RSA}}(\sigma_i||\text{ID}_{M})$ of length $|\sigma_i^{\text{RSA}}|=3072$ bits. 

The SAFE protocol also defines the refresh time as
\begin{equation}
    \Delta^{\text{SAFE}}(M,m,t) = a(M,m,t) + b,
\end{equation}
where $a(M,m,t)$ is an empirically measured coefficient depending on the system parameters, and $b\approx 60 \;\mu s$ accounts for the constant signature verification overhead in SAFE with DS applying RSA-3072 bits. Table~\ref{tab:SAFE_time} (adapted from~\cite{Bouchmann2020SAFE} Fig. 2, part (a) to (c) for the growth of $a(\cdot)$ over the size of the data) shows measured refresh times for different data sizes $|M|$ and $(m,t)$ thresholds. Note that SAFE implementation in \cite{Bouchmann2020SAFE} is coded by JAVA and run over Intel i7-7700 CPU clocked at 3.6 GHz, with 8GB of RAM and Ubuntu 16.04.5 OS.

Let $\Delta = \max \lbrace \frac{\alpha_\text{net}\bigl(|M|+\gamma_{\text{pad}}\bigr)}{k_\text{min}},\;\Delta^\text{SAFE}\rbrace$, where $k_\text{min}$ denotes the minimum required per-link SKR of QKD-based ITS channels, $\alpha_{\text{net}}$ is the network overhead factor that accounts for the number of forward/backward transmissions per $\Delta^\text{Refresh}$ per (TEE, $S_i$) pair and, $\gamma_{\text{pad}}=|\sigma_i^{\text{RSA}}|= 3072\; \forall\; i$ knowing that SAFE runs RSA-based DS and TEE coordinates resharing and signing. So, the maximum number of storable shares of length $|M|$ bits (plus metadata $\gamma_{\text{pad}}$ bits for signatures) is
\begin{equation}\label{eq:delta_refresh}
    N_\text{stor-share}=\Biggl\lfloor \frac{\Delta^{\text{Refresh}}\cdot (365.25\cdot 24 \cdot 3600)}{\Delta} \Biggr\rfloor ,
\end{equation}
 that results in a total transmission budget over one refresh interval equal to
\begin{equation}\label{eq:num_store_share}
    B = N_\text{stor-share}\cdot |M|.
\end{equation}

\begin{table*}[!htb]
     \centering
        \caption{Processing time $a(M,m,t)$ (in seconds) for SAFE operations ``$\mathsf{Share}$''+``$\mathsf{Renew}$''+``$\mathsf{Reconstruction}$'' under different $(M,m,t)$}
     \begin{tabular}{|>{\bfseries}l|*{12}{>{\centering\arraybackslash}p{3.2em}|}}

     \hline
     \diagbox[width=6.em,height=4.5ex]{$(m,t)$}{$|M|$}
 & 100 B & 300 B & 1 KB & 3 KB & 10 KB & 30 KB & 100 KB & 300 KB & 1 MB & 3 MB & 10 MB & 30 MB \\
     \hline\hline
     $(3,2)$ & 14 ms & 18 ms& 19 ms& 21 ms& 31 ms& 77 ms& 240 ms& 650 ms& 2.1 s& 7.5 s& 21 s& 72 s\\
     \hline
     $(5,3)$ & 21 ms & 25 ms& 29 ms& 36 ms& 57 ms& 140 ms& 460 ms& 1.2 s& 4 s& 12 s& 42 s& 130 s \\
     \hline
     $\mathbf{(7,4)}$ & \textbf{30 ms} & \textbf{39 ms}& \textbf{41 ms}& \textbf{49 ms}& \textbf{90 ms}& \textbf{200 ms}& \textbf{700 ms}& \textbf{2 s}& \textbf{5.7 s}& \textbf{18 s}& \textbf{63 s}& \textbf{200 s}\\
     \hline
     $\mathbf{(9,5)}$ & \textbf{41 ms}& \textbf{52 ms}& \textbf{54 ms}& \textbf{68 ms}& \textbf{110 ms}& \textbf{270 ms}& \textbf{830 ms}& \textbf{2.6 s}& \textbf{7.5 s}& \textbf{25 s}& \textbf{82 s}& \textbf{260 s}\\
    \hline
     $(11,6)$& 61 ms & 62 ms& 64 ms& 82 ms& 140 ms& 320 ms& 1 s& 2.9 s& 9 s& 27 s& 105 s& 320 s\\
    \hline
    \end{tabular}
    \label{tab:SAFE_time}
\end{table*}
\subsection{Event Definitions and Risk Probability}\label{subsec:RiskProb}

We define two events: $\mathcal{R}_\text{round}$, the recovery of the secret within a refresh interval $\Delta^{\text{Refresh}}$, and $\mathcal{R}_\text{total}$, the recovery over the total time horizon $T$. Let $p_\text{comp}$ denote the probability of compromise per node (and independent of the other nodes) in one interval $\Delta^{\text{Refresh}}$. Then the probability of secret recovery in a single round is
\begin{equation}\label{eq:p_round}
    \Pr[\mathcal{R}_\text{round}] = \sum_{k=t}^m \binom{m}{k} \, (p_\text{comp})^k \, (1-p_\text{comp})^{m-k}.
\end{equation}
Across the full horizon $T$, the number of rounds is $N_\text{rounds}=\lfloor T/\Delta^{\text{Refresh}}\rfloor$, resulting in the following.
\begin{equation}\label{eq:p_total}
    \Pr[\mathcal{R}_\text{total}] = 1 - \left(1-\Pr[\mathcal{R}_\text{round}]\right)^{N_\text{rounds}}.
\end{equation}
Thus, proactive refresh effectively decouples compromises across rounds, sharply reducing long-term risk.
\subsection{Estimating $p_\text{comp}$}

The parameter $p_\text{comp}$ depends on software/hardware vulnerabilities, insider risks, and physical threats. Conservative values (e.g., $p_\text{comp}=0.01$ per 10 years for nodes based on software guard extensions (SGX)) are often used. Publicly documented attacks (e.g., speculative execution, controlled channel, sealing key extraction) demonstrate that SGX's non-negligible compromise risk\cite{Fei2021SGX}. SAFE enhancements, such as attestation, logging, and hardened TEEs may reduce this to $p_\text{comp}\approx0.001$ for $\Delta^{\text{Refresh}}=10$ years. In practice, empirical estimates may be obtained from security reports (e.g., ENISA, IBM X-Force), mapping annual compromise probabilities into intervals $\Delta^{\text{Refresh}}$. A more data-driven estimate of $p_\text{comp}$ can be derived from the
\emph{Exploit Prediction Scoring System (EPSS)} model \cite{EPSS}, which provides, short-term empirical probabilities of exploitation for individual software/hardware vulnerabilities.
\section{ILP Formalism}
\label{sec:ilp}
We formulate SAFE deployment as a scenario--based ILP with stochastic link costs. 
Let $\Psi$ cost realizations be sampled, indexed by $\psi\in\{1,\ldots,\Psi\}$. 
Each scenario is solved independently and aggregated to obtain feasibility and expected cost.

\paragraph{Sets and parameters.}
Let $\mathcal{N}=\{1,\ldots,n\}$ be candidate sites (each can host a storage node or a hub/TEE). 
We select $m\le n$ storage nodes and $h$ hubs (typically $h=1$). 
Let $d_{i,j}$ be Euclidean distance, $k_{\text{min}}$ the required SKR per storage–hub pair, and 
$[\underline{C^{(\cdot)}_{i,j}},\overline{C^{(\cdot)}_{i,j}}]$ and $k^{(\cdot)}_{i,j}$ denote cost interval and SKR of modality 
$\cdot\!\in\!\{\mathrm{CV},\mathrm{CV{+}MUX},\mathrm{DV},\mathrm{DV{+}MUX}\}$. 
The financial budget is $B_{\text{fin}}$ and $\mathrm{W}=1000\,B_{\text{fin}}$ penalizes violations. The parameter $\mathrm{M}$ denotes a sufficiently large upper bound on the link multiplicity.

\paragraph{Decision variables (per scenario $\psi$).}
The Binary variables $x_i$ and $h_i$ indicate the storage and the hub selection; 
$\mathsf{leaf}_i$ indicates an active storage leaf ($x_i=1,h_i=0$); 
$\mathsf{assign}_{j,i}$ indicates that the storage $j$ is attached to the hub $i$. 
The multiplicities of links are 
$l^{(\cdot)}_{j,i}\!\in\!\mathbb{Z}_{\ge0}$ and classical links $l^{\mathrm{cl}}_{j,i}\!\in\!\{0,1\}$. 
Slack $s^{\psi}\!\ge\!0$ allows for occasional budget violations and $\mathsf{M}\in\mathbb{N}$ is an arbitrary large value.

\paragraph{Cost model.}
Per scenario $\psi$,
\[
\begin{aligned}
\mathsf{Cost}(\psi)&=\!\!\sum_{i\neq j}\!
\Big(C^{\mathrm{CV},\psi}_{i,j}l^{\mathrm{CV}}_{j,i}
+C^{\mathrm{CV{+}MUX},\psi}_{i,j}l^{\mathrm{CV{+}MUX}}_{j,i}\\
&+C^{\mathrm{DV},\psi}_{i,j}l^{\mathrm{DV}}_{j,i}
+C^{\mathrm{DV{+}MUX},\psi}_{i,j}l^{\mathrm{DV{+}MUX}}_{j,i}
+C^{\mathrm{cl},\psi}_{i,j}l^{\mathrm{cl}}_{j,i}\Big),
\end{aligned}
\]
and the objective is $\min\;\mathsf{Cost}(\psi)+\mathrm{W}\,s^{\psi}$.

\paragraph{Constraints.}
\begin{itemize}
\item{\bf Selection:} $\sum_i x_i=m,\;\sum_i h_i=h.$ 
\item{\bf Leaf definition:} $\mathsf{leaf}_i\le x_i,\;\mathsf{leaf}_i\le1-h_i.$ 
\item{\bf Star attachment:} $\sum_i\mathsf{assign}_{j,i}=\mathsf{leaf}_j,$ 
$\mathsf{assign}_{j,i}\le h_i,\;\mathsf{assign}_{j,i}\le x_j.$ 
\item{\bf Link activation:} $l^{(\cdot)}_{j,i}\le \mathsf{M}\cdot\mathsf{assign}_{j,i}$ and 
$l^{\mathrm{cl}}_{j,i}\le\mathsf{assign}_{j,i}.$ 
\item{\bf Capacity:} for any active pair,
\[\begin{aligned}
k_{\text{min}}\le
\Big(&k^{\mathrm{CV}}_{j,i}l^{\mathrm{CV}}_{j,i}
+k^{\mathrm{CV{+}MUX}}_{j,i}l^{\mathrm{CV{+}MUX}}_{j,i}+\\
&k^{\mathrm{DV}}_{j,i}l^{\mathrm{DV}}_{j,i}
+k^{\mathrm{DV{+}MUX}}_{j,i}l^{\mathrm{DV{+}MUX}}_{j,i}\Big).
\end{aligned}\]
\item{\bf Soft budget:} $\mathsf{Cost}(\psi)\le B_{\text{fin}}+s^{\psi}.$
\end{itemize}
\paragraph{Scenario aggregation.}
We report the feasibility rate
\[
\rho=\frac{1}{\Psi}\sum_{\psi=1}^{\Psi}\mathbf{1}\{s^{\psi}=0\},
\]
and average feasible cost $\overline{\mathsf{Cost}^*}$. 
This formulation avoids explicit bounds on link multiplicities; feasibility is enforced through SKR and budget constraints, while slack allows rare violations with heavy penalty.
\section{Integrated Risk-Optimization Framework}\label{sec: integ}
\begin{figure*}[htb!]
\centering
  \includegraphics[width=\textwidth]{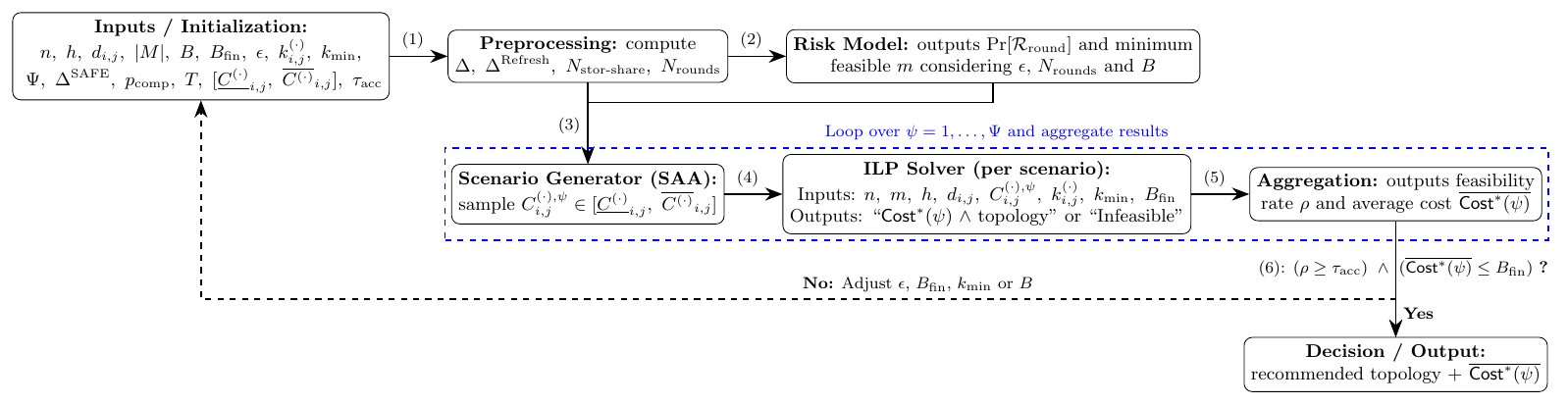} 
  \caption{Diagram of integrated risk and stochastic ILP optimization workflow}
    \label{fig:integ}
\end{figure*}

The flowchart in Fig.~\ref{fig:integ} and the pseudocode in Algorithm~\ref{alg:integrated} jointly summarize the proposed integrated framework, where the risk model provides the minimal number of active storage nodes \(m\) consistent with the security requirement \(\epsilon\), and the stochastic ILP explores multiple cost realizations under bounded uncertainty. The acceptance threshold \(\tau_{\text{acc}}\in[0.8,0.95]\) is introduced as a design parameter that controls the required feasibility rate \(\rho\), i.e., the fraction of scenarios solved within budget. In practice, values closer to \(0.9\) or \(0.95\) provide a robust trade-off between conservatism and deployability. Finally, among the feasible scenarios, the recommended topology is identified as the configuration that occurs most frequently, which ensures consistency between stochastic instances while respecting both budget and risk constraints.
\begin{algorithm}[htb!]
\caption{Integrated Risk--Optimization Pipeline}
\label{alg:integrated}
\DontPrintSemicolon
\KwIn{$n,h,d_{i,j},|M|,B,B_{\text{fin}},\epsilon,k^{(\cdot)}_{i,j},k_{\text{min}},\Psi,\Delta^{\text{SAFE}},p_{\text{comp}},$\\
\hspace{1.5em}$T,[\underline{C^{(\cdot)}},\overline{C^{(\cdot)}}]$}
\KwOut{recommended topology, $\rho$, $\overline{\mathsf{Cost}^*}$}
Preprocessing: compute $\Delta,N_{\text{stor-share}}\ (\text{eq. \ref{eq:num_store_share}}),$\; 
$\Delta^{\text{Refresh}}\ (\text{eq. \ref{eq:delta_refresh}}),N_{\text{rounds}}$\;
Risk: invert $\epsilon$ to obtain $\Pr[\mathcal{R}_\text{round}]\ (\text{eq. \ref{eq:p_total}})$ and determine feasible $m\ (\text{eq. \ref{eq:p_round}})$\;
\For{$\psi\leftarrow 1$ \KwTo $\Psi$}{
  Sample $C^{(\cdot),\psi}_{i,j}\in[\underline{C^{(\cdot)}}_{i,j},\overline{C^{(\cdot)}}_{i,j}]$ uniformly at random\;
  Build scenario ILP with variables $l^{(\cdot)}_{j,i},l^{\text{cl}}_{j,i},\mathrm{s}^{\psi}$\;
  Add capacity, assignment and classical-link constraints\;
  Define $\mathsf{Cost}(\psi)$ as linear form\;
  Add soft-budget: $\mathsf{Cost}(\psi)\le B_{\text{fin}}+\mathrm{s}^{\psi}$\;
  Solve $\min\ \mathsf{Cost}(\psi)+\mathrm{W}\,\mathrm{s}^{\psi}$ and record feasible/infeasible solution\;
}
Compute $\rho\leftarrow\frac{1}{\Psi}\sum_{\psi}\mathbf{1}\{\mathrm{s}^{\psi}=0\}$\;
Compute $\overline{\mathsf{Cost}^*}$ over feasible scenarios\;
Decision: if $\rho\ge\tau_{\text{acc}}$ and $\overline{\mathsf{Cost}^*}\le B_{\text{fin}}$ then accept else advise relaxation\;
\end{algorithm}

\section{Numerical Results}
\label{sec:synth}

This section evaluates the proposed SAFE optimization framework through two complementary studies. First, randomized synthetic networks are considered to assess the influence of the deployment scale, transmission budget, and QKD parameters under statistically generated topologies. Second, the framework is validated on realistic metropolitan-scale QKD infrastructures derived from the Greater Paris area, demonstrating its applicability to practical deployments.

\subsection{Randomized Synthetic Network Evaluation}

\subsubsection{Simulation Setup}

The first set of experiments aims to characterize the behavior of the proposed optimization framework under generic deployment conditions. Candidate storage-node locations are generated uniformly at random within square fields of side length $X$, allowing the impact of the geographical scale to be evaluated independently of any particular network topology.

An ideal transmission regime following (\ref{eq:capacity}) is assumed throughout this study. The number of candidate storage nodes is fixed to either $n=12$ or $20$, with
$X\in\{60,140,180,220,300,380,700\}$ km. For each field realization, the integrated optimization framework of Algorithm~\ref{alg:integrated} is executed over $N_{\mathrm{TRIALS}}=300$ independent random deployments and $\Psi=20$ scenario realizations obtained by sampling the modality-dependent deployment costs uniformly within their corresponding intervals.

Unless otherwise stated, the financial budget is fixed to
$B_{\mathrm{fin}}=10^{10}\,\$$, the node compromise probability to
$p_{\mathrm{comp}}=10^{-4}$ per year, the storage lifetime to
$T=100$ years, and the penalty multiplier to
$\mathrm{W}=100$. The transmission budget $B$, $k_{\min}$, the global compromise tolerance $\epsilon$, and the resulting optimal number of active storage nodes $m$ are specified in each figure. The value of $m$ is automatically determined by the integrated risk model presented in Section~\ref{subsec:RiskProb}.

To investigate the influence of the relative deployment costs of DV- and CV-QKD technologies, three representative cost regimes are considered throughout the randomized-network evaluation.

\begin{itemize}
    \item \textbf{Small-gap regime:} CV $\gamma=[1.0,1.3]$, DV $\gamma=[0.7,1.0]$, and classical $\gamma=[0.95,1.05]$.
    
    \item \textbf{Large-gap regime:} CV $\gamma=[0.7,1.0]$, DV $\gamma=[1.0,1.3]$, and classical $\gamma=[0.95,1.05]$.
    
    \item \textbf{No-gap regime:} CV $\gamma=[0.95,1.05]$, DV $\gamma=[0.95,1.05]$, and classical $\gamma=[0.95,1.05]$.
\end{itemize}

Unless otherwise specified, the physical-layer capacities follow (\ref{eq:capacity}), and the considered cost regime is explicitly indicated in the corresponding figure.

\subsubsection{Results and Key Observations}

\begin{figure}[t]
\centering
\includegraphics[height=0.3\textheight,width=0.49\textwidth]{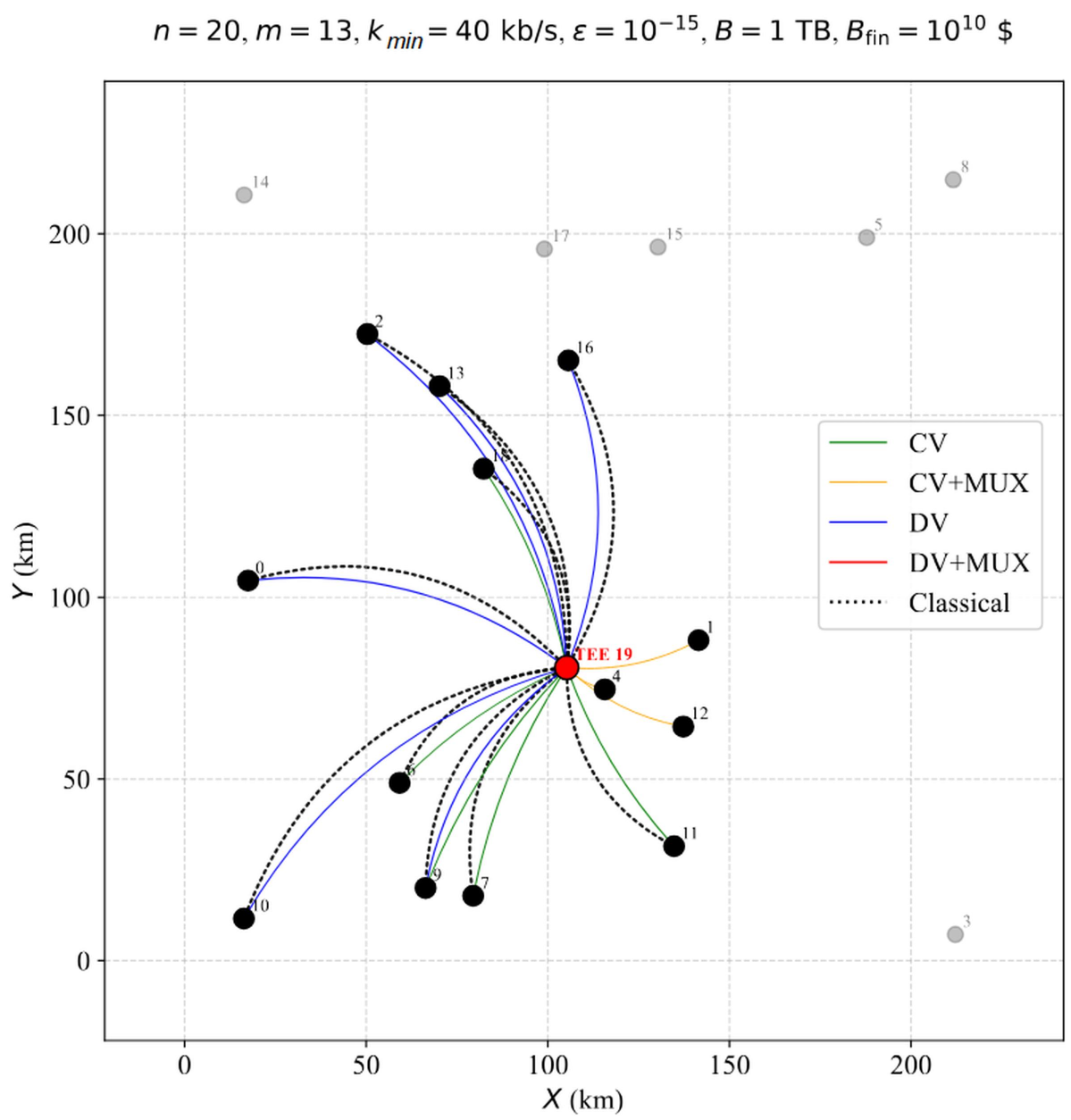}
\caption{An optimized SAFE topology sample under no-gap regime and inside a squared field of side $X=220$ km.}
\label{fig:optimized_topology}
\end{figure}

\begin{figure*}[htb!]
\centering
\begin{tabular}{c@{~~}c}
\includegraphics[width=0.45\textwidth]{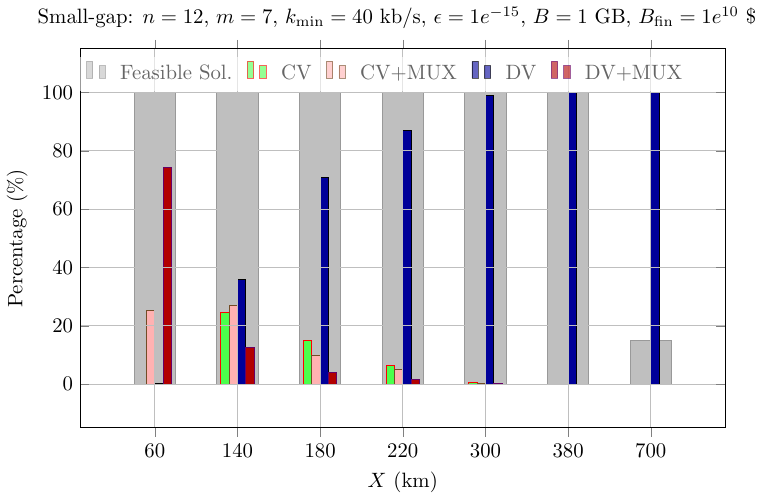} &
\includegraphics[width=0.45\textwidth]{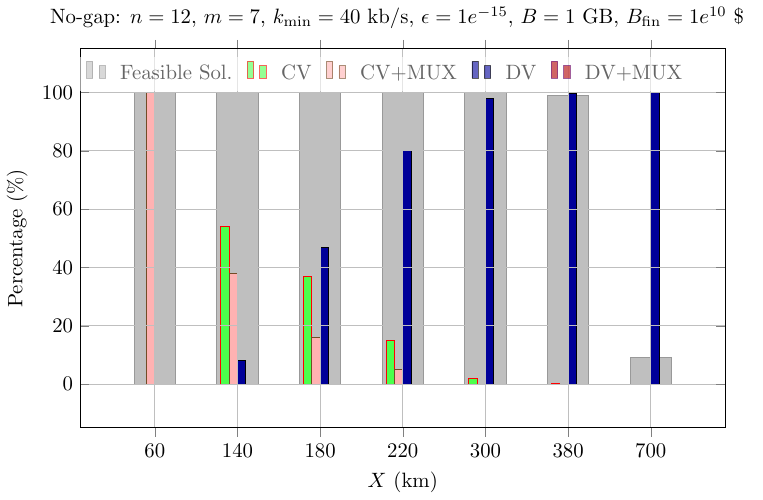} \\ 
~~(a) & ~~(b) \\ ~&~\\
\includegraphics[width=0.45\textwidth]{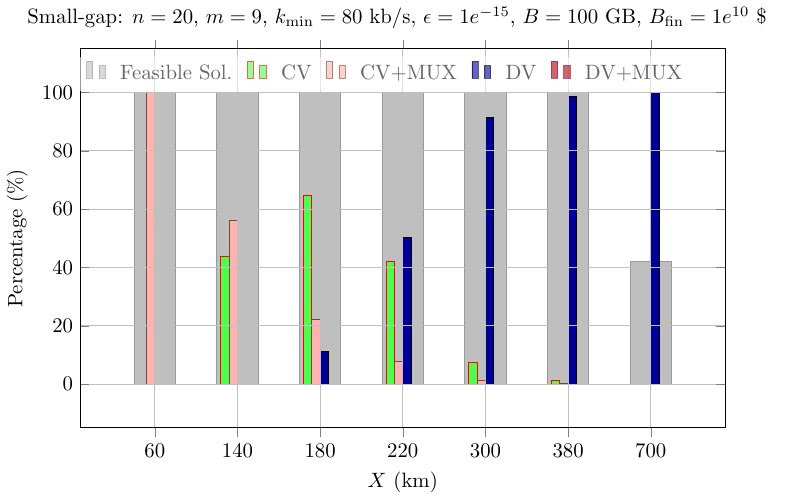} &
\includegraphics[width=0.45\textwidth]{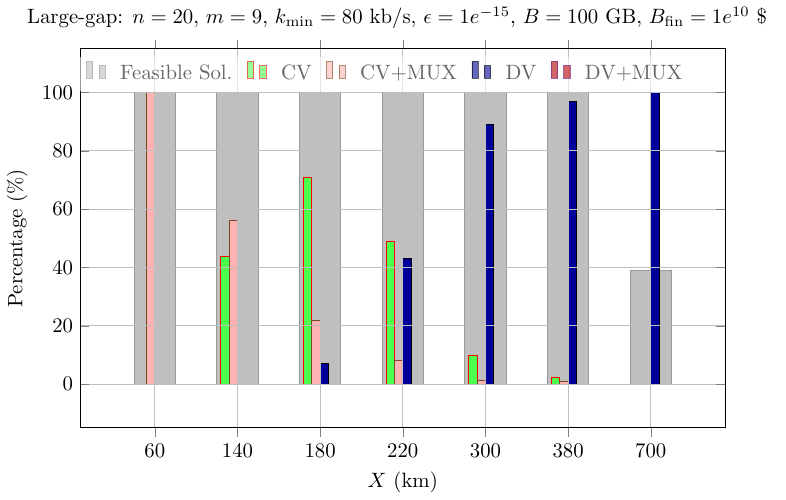} \\ 
~~(c) & ~~(d) \\ ~&~\\
\includegraphics[width=0.45\textwidth]{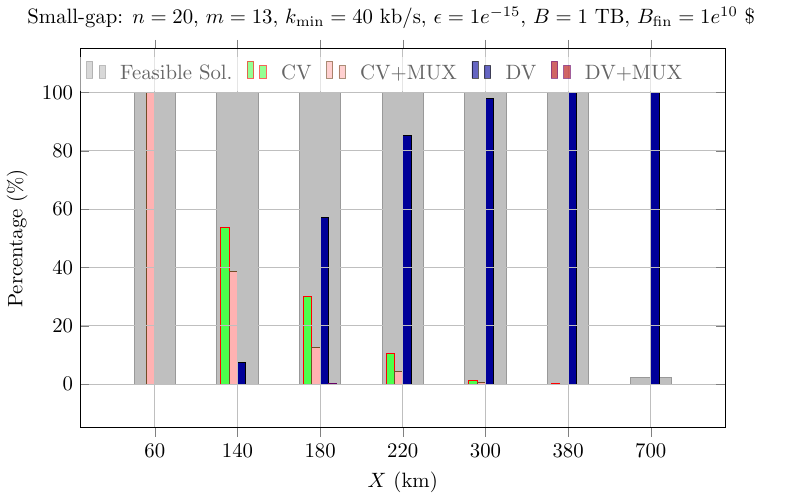} &
\includegraphics[width=0.45\textwidth]{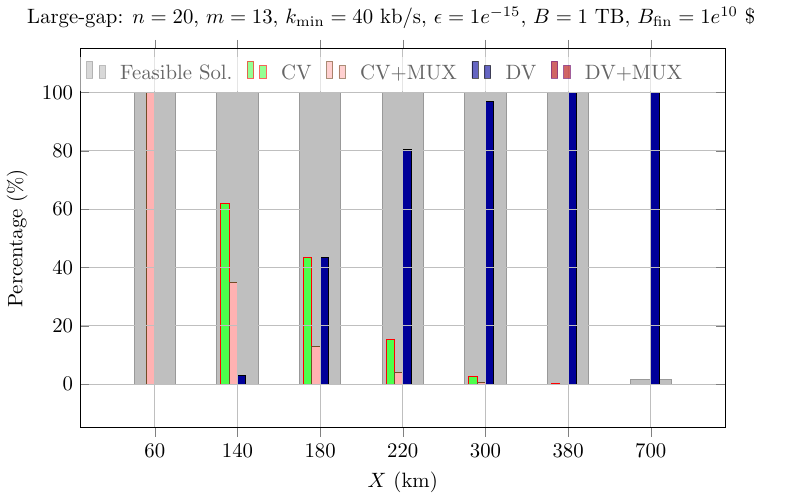} \\ 
~~(e) & ~~(f)
\end{tabular}
\caption{Field-side length ($X$) in km versus the percentage of modality usage CV, DV, CV$+$MUX, DV$+$MUX and the feasibility rate $\rho$ (thick-gray bar). Each case corresponds to distinct sets of $(m, k_{\text{min}}, \epsilon, B, B_{\text{fin}})$ as indicated on top of each subfigure.}
\label{fig:percentage}
\end{figure*}

\begin{figure*}[htb!]
\centering
\begin{tabular}{c@{~~}c}
\includegraphics[width=8.9cm, height=7cm]{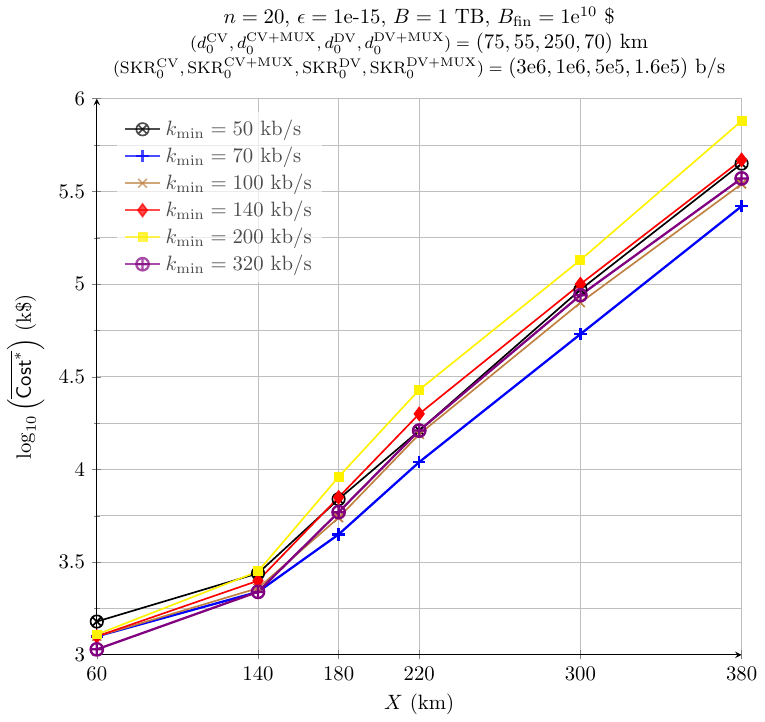} &
\includegraphics[width=8.9cm, height=7cm]{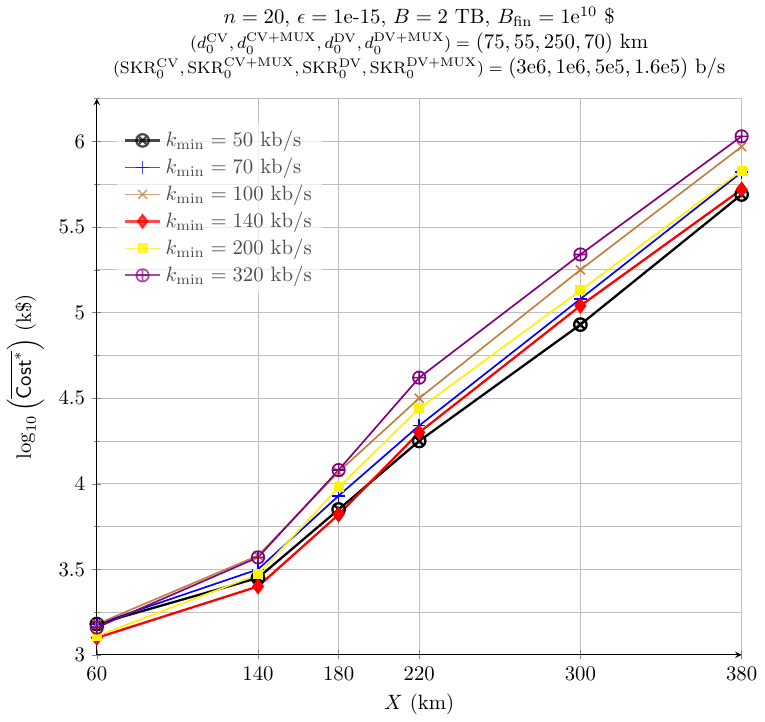} \\ 
~~(a) & ~~(b) \\~&~\\
\includegraphics[width=8.9cm, height=7cm]{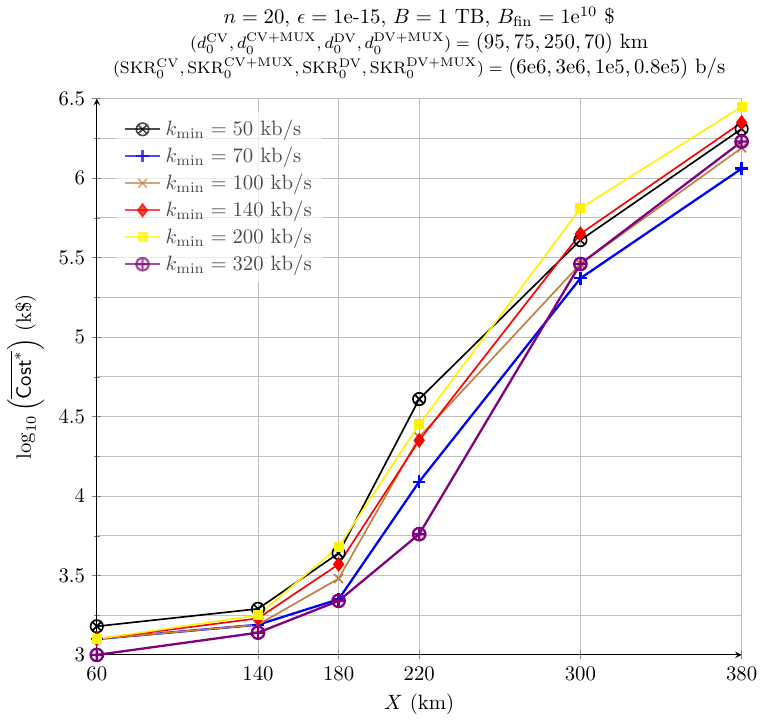} &
\includegraphics[width=8.9cm, height=7cm]{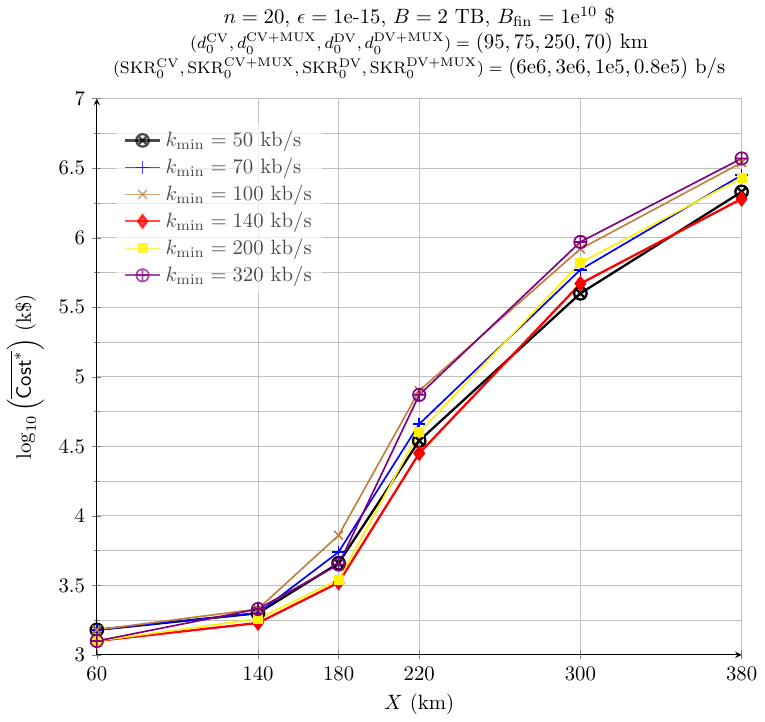} \\ 
~~(c) & ~~(d)
\end{tabular}
\caption{Average deployment costs ($\log_{10}$ scale) versus field-side length $X$ for different minimum required per-link SKRs $k_\text{min}$ are plotted. Only the all-modalities scenario ``CV $\&$ DV $\&$ CV$+$MUX $\&$ DV$+$MUX'' in ``no-gap regime'' is considered. SKRs at zero distance as well as the cut-off distances are reported on top of each sub-figure.}
\label{fig:Descrete_Average_Cost}
\end{figure*}

The results of the randomized synthetic-network study are summarized in Figs.~\ref{fig:optimized_topology}--\ref{fig:Descrete_Average_Cost}. Together, these figures illustrate how the proposed SAFE optimization jointly determines the network topology, the preferred QKD modality, and the resulting deployment cost under different geographical scales and security requirements.

Figure~\ref{fig:optimized_topology} presents a representative optimal SAFE deployment obtained under the no-gap regime for a square field of side length $X=220$ km. The synthesized topology illustrates the spatial distribution of the selected active storage nodes and the corresponding QKD links. Depending on the transmission distance and the associated deployment costs, the optimization naturally combines different QKD modalities to minimize the total infrastructure cost while satisfying the required security constraints.

The influence of the geographical scale and the deployment-cost model is further investigated in Fig.~\ref{fig:percentage}. The gray bars indicate the feasibility rate $\rho$, while the colored curves report the percentage of links assigned to each QKD modality. As the field size $X$ increases, longer transmission distances progressively reduce the feasible operating region and modify the optimal modality allocation. The transition between CV, DV and multiplexed solutions therefore depends not only on the transmission distance but also on the assumed relative deployment costs, the transmission budget $B$, and the required minimum per-link secret-key rate $k_{\min}$.

Finally, Fig.~\ref{fig:Descrete_Average_Cost} analyzes the sensitivity of the optimal deployment cost to the geographical scale, the transmission budget, the minimum required SKR, and the physical-layer parameters $(\mathrm{SKR}_0^{(\cdot)},d_0^{(\cdot)})$. The reported costs correspond to the average optimal ILP solution over $300$ independent network realizations under the no-gap regime, considering the coexistence of all four QKD modalities.

The main observations can be summarized as follows.

\begin{itemize}
\item The average deployment cost increases with the deployment scale $X$, reflecting the longer QKD links and the additional infrastructure required to maintain secure connectivity.

\item Increasing the transmission budget $B$ consistently increases the infrastructure cost because larger amounts of secret material require higher key-generation throughput throughout the network.

\item The preferred QKD modality is jointly determined by the geographical scale, the relative deployment costs, and the physical-layer capabilities, resulting in distinct transition regions where different technologies become economically preferable.

\item Improving the physical-layer parameters $(\mathrm{SKR}_0^{(\cdot)},d_0^{(\cdot)})$ significantly mitigates the infrastructure cost by extending the feasible transmission range of QKD links.

\item \textit{Most importantly, increasing the minimum required per-link secret-key rate $k_{\min}$ does not necessarily reduce the deployment cost.} Once the optimal modality allocation has been reached, the infrastructure cost is primarily governed by the combined interaction between the network topology, the selected QKD technologies, and their transmission reach, rather than by the achievable SKR alone.
\end{itemize}

\subsection{Paris Metropolitan Network Evaluation}

\subsubsection{Simulation Setup}

The second set of experiments validates the proposed SAFE optimization framework on realistic metropolitan-scale QKD infrastructures representative of the Greater Paris area. Unlike the randomized study, the network topology is fixed and corresponds to actual academic laboratories, industrial facilities, and data-center locations, enabling the evaluation of practical deployment scenarios.

Three representative network configurations are considered: a Paris Metro-Scale network with $n=12$ candidate storage sites, a Sparse Greater Paris network with $n=15$ sites, and a Dense Greater Paris network with $n=25$ sites. The names and geographical locations of all candidate sites are listed in Appendix~\ref{app:network_description}. Fiber-path distances are approximated from the Euclidean distances using a detour factor $\alpha=2.0$ to account for realistic fiber routing.

The industrial implementation and multiplexed coexistence model described by (\ref{IndustSKR}) is adopted throughout this study. Since both the node coordinates and the physical link distances are fixed, the integrated optimization framework of Algorithm~\ref{alg:integrated} is executed for a single deterministic scenario ($\Psi=1$), eliminating the statistical averaging required in the randomized evaluation. Consequently, each experiment produces a unique optimal deployment cost together with the corresponding SAFE topology. The parameters $B$, $k_{\min}$, $\epsilon$, and  $m$ are reported in each figure.

\begin{figure*}[htb!]
\centering
\begin{tabular}{c@{~~}c}
\includegraphics[width=8.5cm, height=7.7cm]{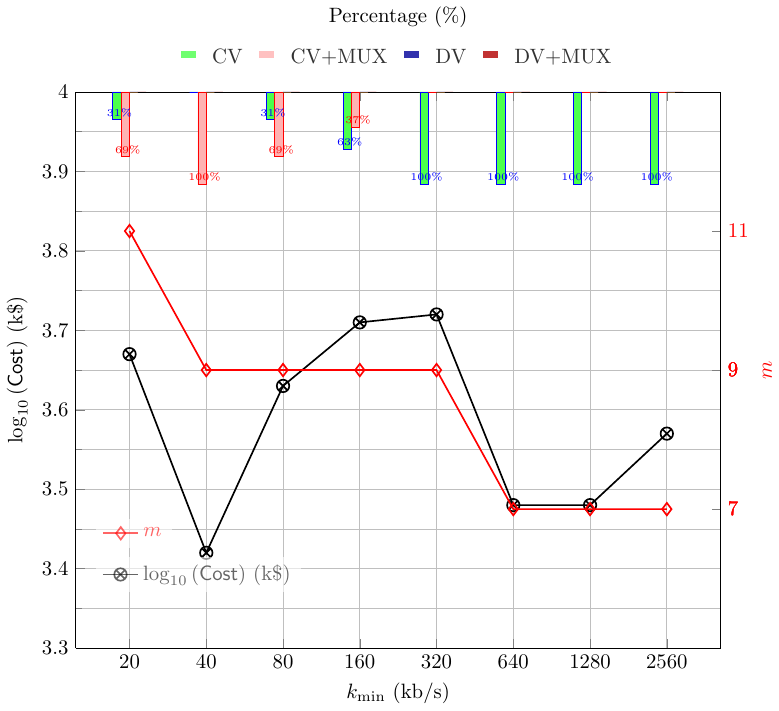} &
\includegraphics[width=8.6cm, height=7.4cm]{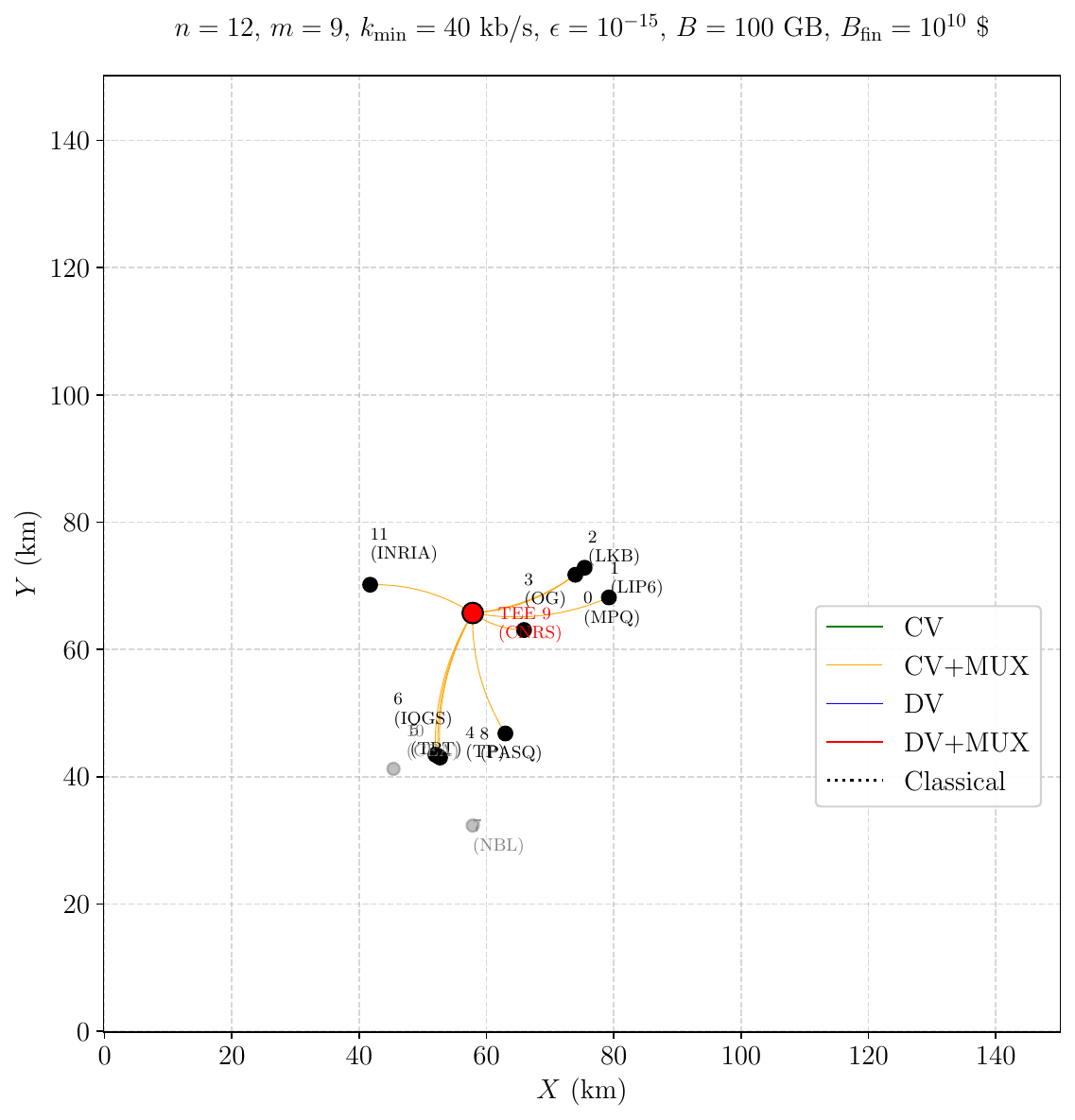} \\ 
~~(a) & ~~(b)
\end{tabular}
    \caption{SAFE optimization for the Paris Metro-Scale network with $12$ candidate storage sites. Subfigure (a) shows the total QKD deployment cost (left axis), the number of selected storage nodes m (right axis), and the percentage of each QKD modality versus the required minimum per-link secret-key rate $k_{min}$. Subfigure (b) depicts the corresponding minimum-cost SAFE topology in (a), where $(x,y)$ denotes coordinate of sites.}
    \label{fig:ParisMetroScaleNet}
\end{figure*}
\begin{figure*}[htb!]
\centering
\begin{tabular}{c@{~~}c}
\includegraphics[width=8.5cm, height=7.7cm]{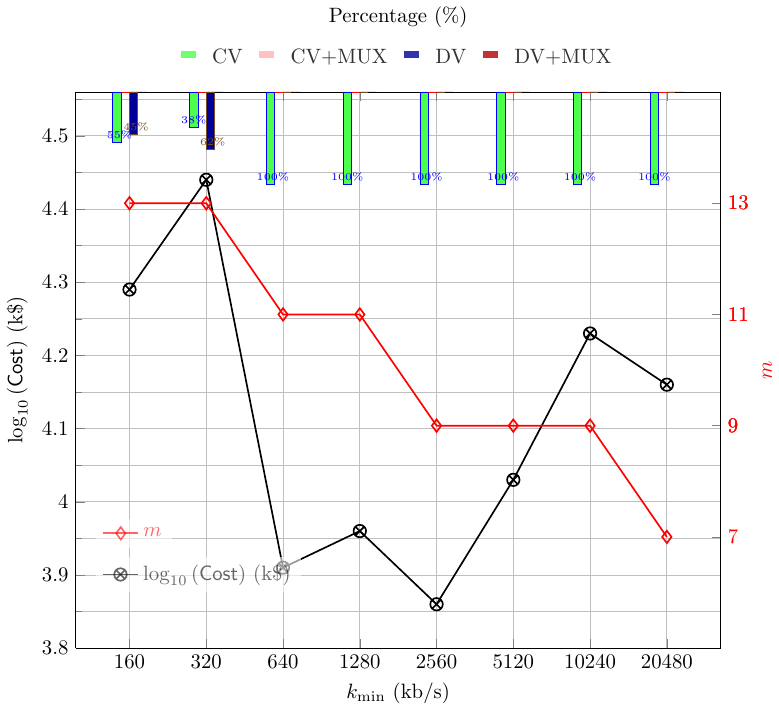} &
\includegraphics[width=8.6cm, height=7.4cm]{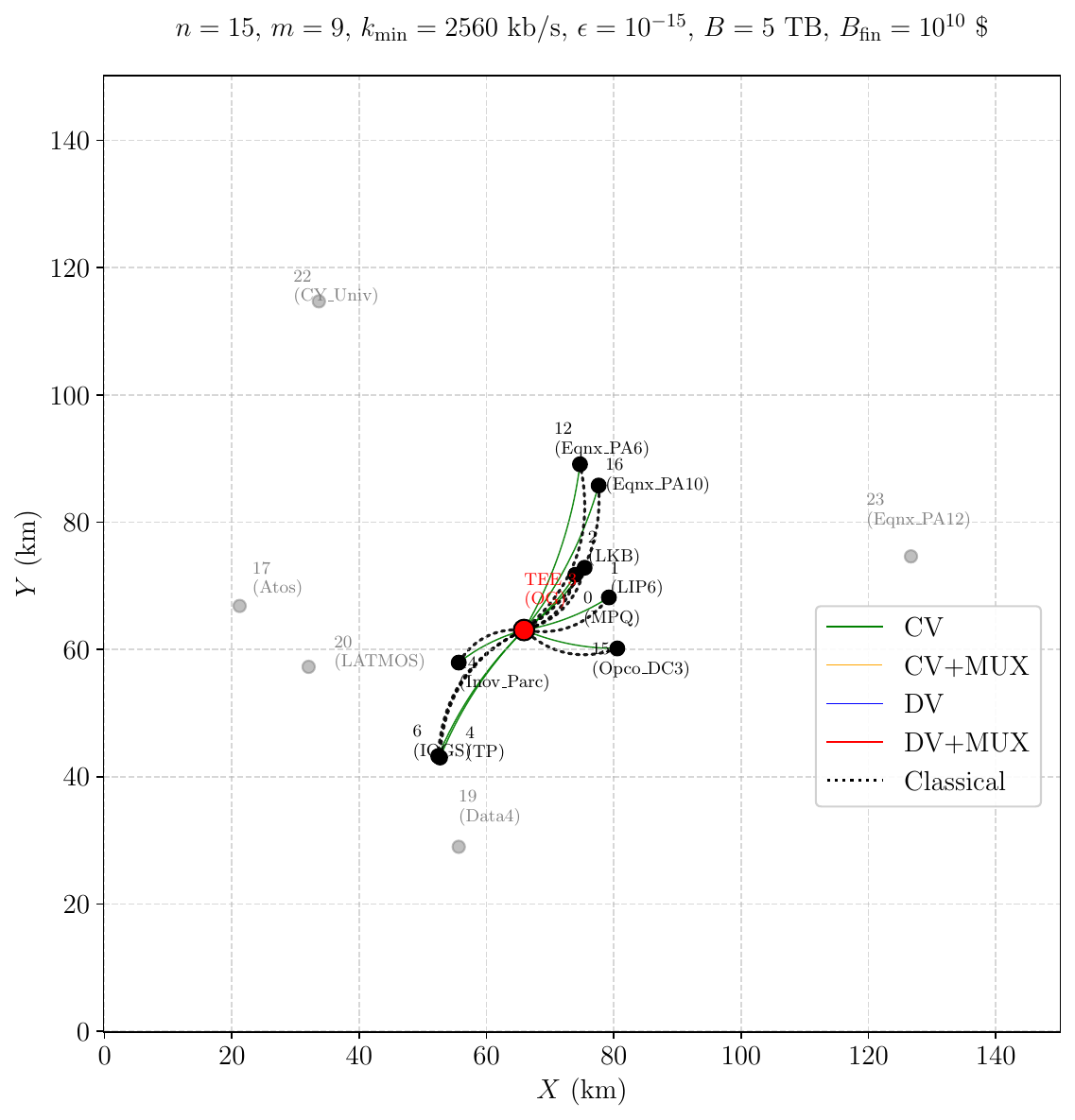} \\ 
~~(a) & ~~(b)
\end{tabular}
    \caption{SAFE optimization for the Sparse Greater Paris network with $15$ candidate storage sites. Subfigure (a) reports the total QKD deployment cost (left axis), the selected number of storage nodes m (right axis), and the QKD modality allocation as functions of the minimum required per-link secret-key rate $k_{min}$. Subfigure (b) depicts the corresponding minimum-cost SAFE topology in (a), where $(x,y)$ denotes coordinate of sites.}
    \label{fig:SparseGreaterParisNet}
\end{figure*}

\subsubsection{Results and Key Observations}

\begin{figure*}[htb!]
\centering
\begin{tabular}{c@{~~}c}
\includegraphics[width=8.5cm, height=7.7cm]{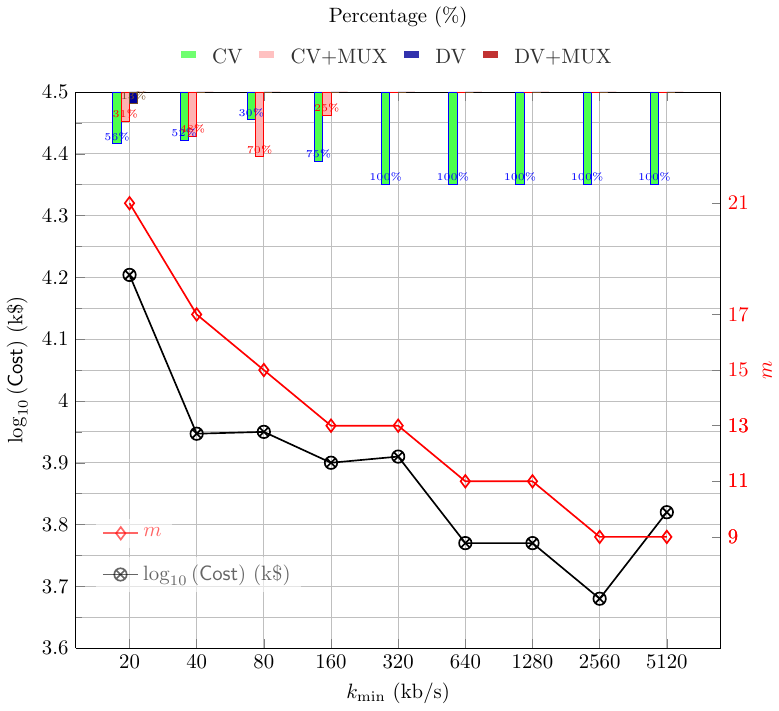} &
\includegraphics[width=8.6cm, height=7.4cm]{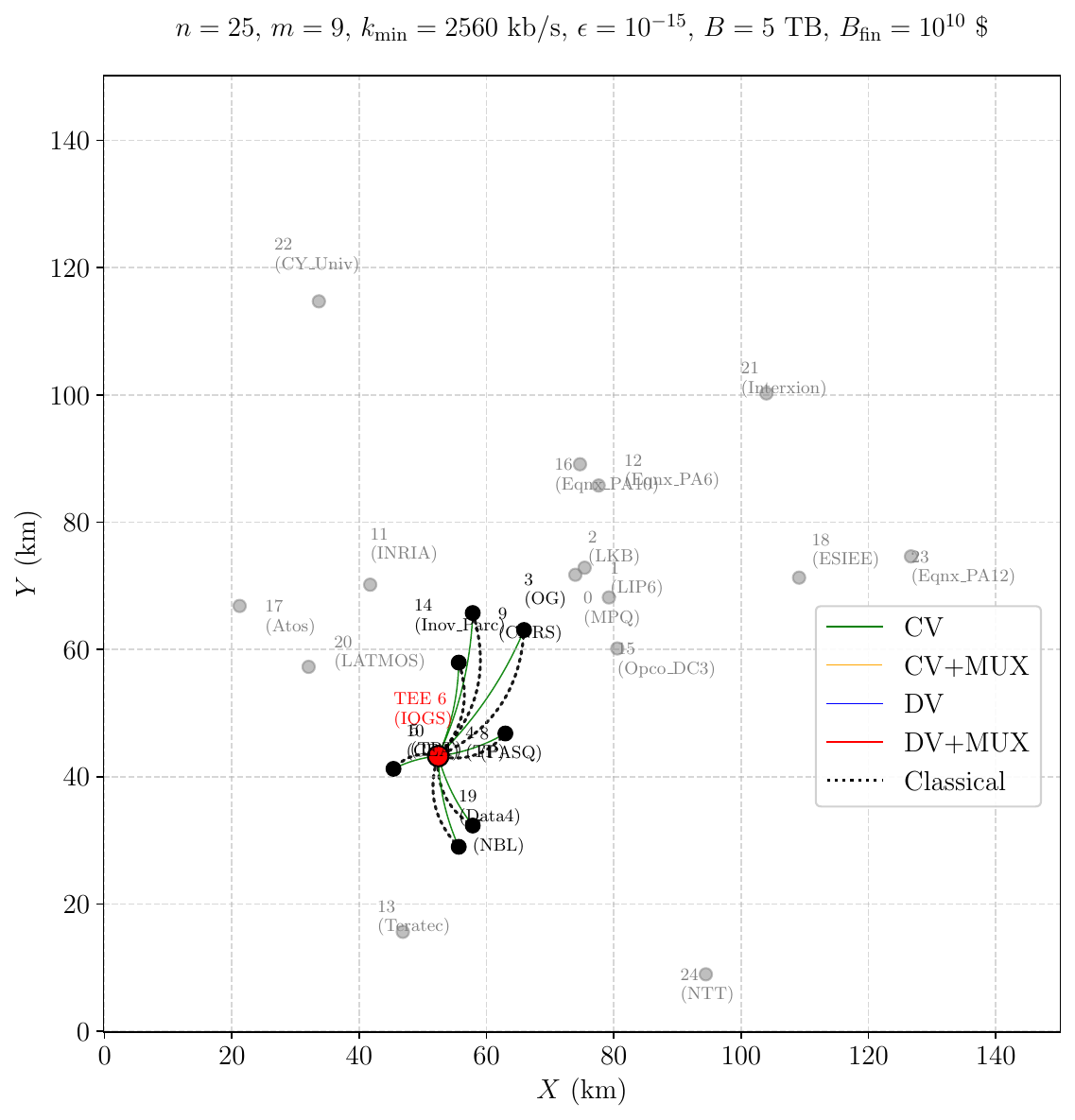} \\ 
~~(a) & ~~(b)
\end{tabular}
    \caption{SAFE optimization for the Dense Greater Paris network with $25$ candidate storage sites. Subfigure (a) illustrates how the total QKD deployment cost, the optimal number of storage nodes m, and the modality allocation evolve with the minimum required per-link secret-key rate $k_{min}$.  Subfigure (b) depicts the corresponding minimum-cost SAFE topology in (a), where $(x,y)$ denotes coordinate of sites.}
    \label{fig:DenseGreaterParisNet}
\end{figure*}

The optimization results for the three metropolitan-scale QKD networks are presented in Figs.~\ref{fig:ParisMetroScaleNet}--\ref{fig:DenseGreaterParisNet}. For each network, subfigure (a) jointly reports the total QKD deployment cost, the resulting number of active storage nodes $m$, and the percentage allocation of the four QKD modalities as functions of the $k_{\min}$. The corresponding minimum-cost SAFE topology is illustrated in subfigure (b) for the highlighted operating point.

Although the three network scenarios differ in terms of geographical coverage, node density, and transmission budget, they exhibit remarkably consistent optimization behavior. The integrated framework jointly determines the optimal number of active storage nodes, the preferred QKD modality, and the resulting network topology in order to minimize the overall deployment cost while satisfying the SAFE security constraints.

The main observations can be summarized as follows.

\begin{itemize}

\item The deployment cost exhibits a \emph{non-monotonic} dependence on the minimum required per-link secret-key rate $k_{\min}$. Contrary to intuition, increasing the achievable SKR does not necessarily reduce the infrastructure cost.

\item For small values of $k_{\min}$, long transmission distances often require the use of longer-reach but more expensive DV-QKD links, or hybrid deployments combining several QKD modalities. As $k_{\min}$ increases, higher-capacity CV-based solutions progressively become sufficient to satisfy the traffic demand, yielding a lower overall deployment cost.

\item Beyond a technology-dependent threshold, further increasing $k_{\min}$ may reduce the number of selected storage nodes $m$, but does not necessarily decrease the total deployment cost. Once the optimal modality allocation has been reached, the deployment cost is primarily governed by the network topology and the transmission reach of the available QKD technologies.

\item The Paris Metro-Scale network reaches its minimum deployment cost around $k_{\min}\approx40$ kb/s using $m=9$ active storage nodes, where the optimal solution relies exclusively on the CV+MUX modality.

\item Both the Sparse and Dense Greater Paris networks achieve their minimum deployment costs around $k_{\min}\approx2560$ kb/s with $m=9$ active storage nodes. Despite their different network sizes, the optimal deployments converge toward predominantly CV-based architectures, highlighting the economic advantage of high-capacity CV-QKD links for metropolitan-scale SAFE infrastructures.

\end{itemize}
\section{Conclusion and Future Work}
\label{sec:concl}
We presented an integrated risk–optimization framework for designing cost-efficient, information-theoretic secure LTSS infrastructures under the SAFE protocol. By coupling a probabilistic compromise model with a stochastic ILP formulation, the method jointly captures long-term confidentiality guarantees and infrastructure cost trade-offs across heterogeneous QKD modalities. Numerical results demonstrate that hybrid CV/DV deployments—particularly with wavelength multiplexing—substantially reduce the average infrastructure cost while maintaining the required global compromise tolerance. The optimized SAFE topologies automatically adapt modality selection to field size, transmission budget, and SKR regimes, confirming that hybrid quantum infrastructures are a practical path toward metropolitan-scale, post-quantum-resilient LTSS systems.

An important direction for future work is to extend and evaluate SAFE toward \emph{dynamic} and \emph{verifiable} secret-sharing settings that arise in real deployments where $n$ is limited and nodes may fail or be compromised before the storage horizon $T$. Dynamic Proactive Secret Sharing (DPSS)\cite{DPSS2015} and Dynamic-Committee PSS (DC-PSS) \cite{DC-PSS2025} enable efficient replacement of compromised nodes and secure migration of shares when committee membership changes, while preserving proactive confidentiality guarantees. Integrating such mechanisms with SAFE would allow rapid reconfiguration of active sets of size $m\leq n$ within intervals shorter than $\Delta^{\text{Refresh}}$, improving robustness against mobile adversaries and denial-of-service events. In parallel, Verifiable Secret Sharing (VSS) techniques—e.g., commitment-based schemes such as PCS \cite{Pedersen1991}—can complement TEE-based authentication by detecting dealer frauds and ensuring share correctness. 

Combining SAFE with DPSS/DC-PSS and VSS therefore represents a promising research avenue toward adaptive, large-scale LTSS infrastructures, and a cost-effectiveness analysis of such solutions is required for the security providers.

\appendices
\section{Dense Greater-Paris Synthetic Quantum-Network}
\label{app:network_description}

List of 25 nodes that correspond to key academic, industrial and data center locations in the Greater Paris area. This extensive but synthetic network is designed to benchmark quantum technologies, specifically QKD in both DV and CV modalities.

\begin{enumerate}[label=\arabic*:, leftmargin=2em, start=0]
    \item \textbf{MPQ (Univ Paris Cité)}: 10 Rue Alice Domon, 75013 Paris.
    \item \textbf{LIP6 (Sorbonne Univ)}: 4 Place Jussieu, 75005 Paris.
    \item \textbf{LKB (Kastler-Brossel)}: 24 Rue Lhomond, 75005 Paris.
    \item \textbf{OG (Orange)}: 44 Ave de la République, 92320 Châtillon.
    \item \textbf{TP (Télécom Paris)}: 19 Place M. Perey, 91120 Palaiseau.
    \item \textbf{TRT (Thales RT)}: 1 Ave Augustin Fresnel, 91767 Palaiseau.
    \item \textbf{IOGS (Inst. Optique)}: 2 Ave Augustin Fresnel, 91127 Palaiseau.
    \item \textbf{NBL (Nokia Bell Labs)}: Route de Villejust, 91620 Nozay.
    \item \textbf{PASQ (Pasqal)}: 7 Rue Léonard de Vinci, 91300 Massy.
    \item \textbf{CNRS}: 1 Place Aristide Briand, 92190 Meudon.
    \item \textbf{CEA}: Bâtiment 466, 91191 Gif-sur-Yvette.
    \item \textbf{INRIA}: Domaine de Voluceau, 78150 Le Chesnay.
    \item \textbf{Eqnx\_PA6}: 10 Rue Waldeck Rochet, 93300 Aubervilliers.
    \item \textbf{Teratec}: 2 Rue de la Piquetterie, 91680 Bruyères-le-Châtel.
    \item \textbf{Inov\_Parc}: Ave de l'Europe, 78140 Vélizy-Villacoublay.
    \item \textbf{Opco\_DC3}: 61 Rue Julian Grimau, 94400 Vitry-sur-Seine.
    \item \textbf{Eqnx\_PA10}: 114 Rue Ambroise Croizat, 93200 Saint-Denis.
    \item \textbf{Atos}: Rue Jean Jaurès, 78340 Les Clayes-sous-Bois.
    \item \textbf{ESIEE}: 2 Blvd Blaise Pascal, 93160 Noisy-le-Grand.
    \item \textbf{Data4}: Route de Nozay, 91460 Marcoussis.
    \item \textbf{LATMOS}: 11 Blvd d'Alembert, 78280 Guyancourt.
    \item \textbf{Interxion}: 1 Rue de l'Haye, 93290 Tremblay.
    \item \textbf{CY\_Univ}: 33 Blvd du Port, 95000 Cergy.
    \item \textbf{Eqnx\_PA12}: 16 Ave J. Froelicher, 77600 Ferrières.
    \item \textbf{NTT (DATA PARIS 1)}: 4 Ave de l'Europe, 91830 Le Coudray-Montceaux.
\end{enumerate}

\section*{Acknowledgment}

...

\bibliographystyle{IEEEtran}
\bibliography{refs}

\begin{IEEEbiography}[{\includegraphics[width=1.1in,height=1.25in,clip]{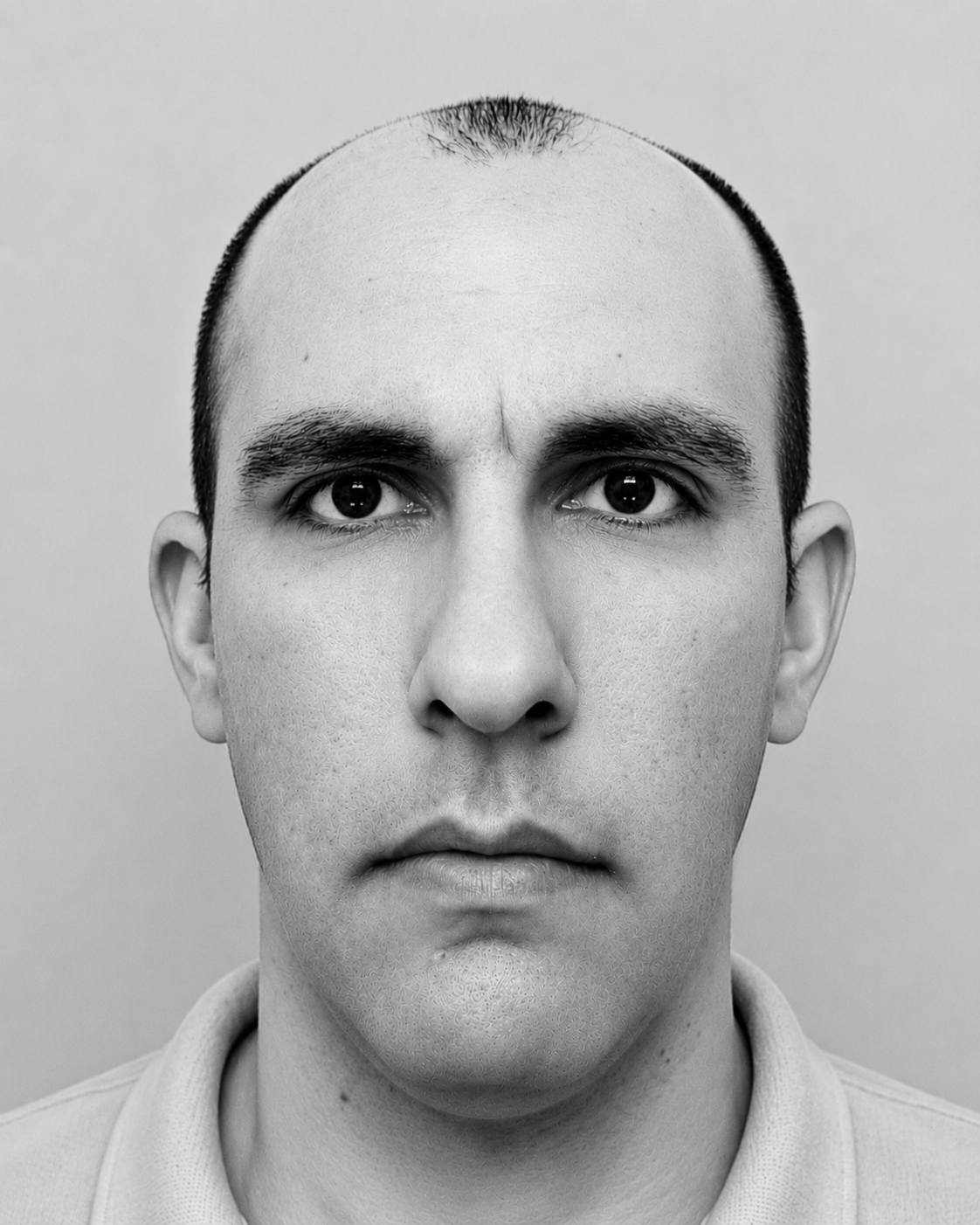}}]{Alireza Tasdighi}
received the Ph.D. degree in applied mathematics, information theory, and coding theory from the Amirkabir University of Technology, Tehran, in 2016. From 2014 to 2015, he was a Visiting Researcher with Carleton University, Ottawa, ON, Canada.

From 2012 to 2019, he was an Assistant Professor with the Institute for Advanced Studies in Basic Sciences (IASBS) and Persian Gulf University, Iran. He subsequently worked as a Post-Doctoral Researcher with LAB-STICC, Universit\'e Bretagne Sud, Lorient, France, from 2019 to 2020, and with IMT Atlantique, Plouzan\'e, France, from 2020 to 2022. From 2022 to 2025, he was a Researcher with the Department of Communication and Electronics (COMELEC), T\'el\'ecom-Paris, Paris, France. He is currently a Scientific Researcher with D\'epartement Informatique et R\'eseaux (INFRES), T\'el\'ecom-Paris. His research interests include information theory, error-correcting codes (LDPC, SC-LDPC), optical fiber communications, machine learning, and quantum and post-quantum cryptography architectures, with a focus on risk modeling and cost optimization for long-term secure quantum key distribution (QKD) storage networks.
\end{IEEEbiography}

\begin{IEEEbiography}[{\includegraphics[width=1.1in,height=1.25in,clip]{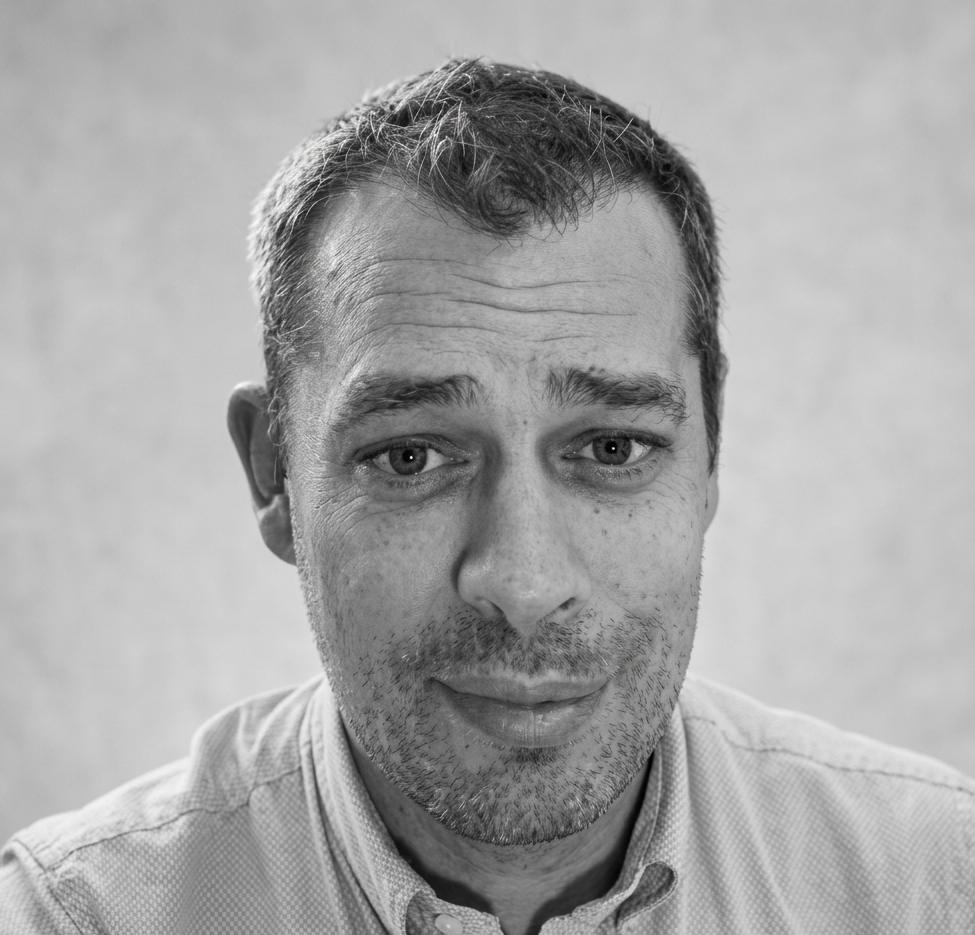}}]{Romain All\'eaume}
graduated from ENS Paris in~2002 and from T\'el\'ecom Paris in~2004. He completed his PhD from Paris VI University in~2004 and was co-recipient of ``magazine La Recherche'' scientific prize. He then joined T\'el\'ecom Paris and coordinated the activity on QKD performed within the European FP6 project SECOQC. Romain All\'eaume co-founded the start-up company SeQureNet that developed, and successfully commercialized the first continuous-variable quantum key distribution (CV-QKD) system in~2013. Currently he is actively involved in the European Quantum Technology Flagship (OpenQKD, Quantum Secure Network Partnership), the EuroQCI program (ParisRegionQCI, FranceQCI) as well as in French PEPR Quantique (QCommTestbed). He is leading a research activity on quantum and classical cryptography, quantum networks and photonic quantum information processing as well as on cross-disciplinary topics linking quantum information with cryptography, photonics and security.
\end{IEEEbiography}

\EOD

\end{document}